\documentclass{aic}

\usepackage{indentfirst}
\usepackage{hyperref}
\usepackage{amssymb, amsmath}
\usepackage{array}
\usepackage{graphicx}
\usepackage[labelsep=none]{caption}

\hypersetup{pdfstartview=FitH, pdfpagemode=UseNone}

\newcounter{question}

\newenvironment{question*}[1]{\emph{Question #1.}}{}
\newcounter{conjecture}
\newenvironment{conjecture}{\refstepcounter{conjecture} \emph{Conjecture \theconjecture.}}{}
\newenvironment{conjecture*}[1]{\emph{Conjecture #1.}}{}
\newcounter{lemma}
\newenvironment{lemma}{\refstepcounter{lemma} \emph{Lemma \thelemma.}}{}
\newenvironment{lemma*}[1]{\emph{Lemma #1.}}{}
\newcounter{proposition}
\newenvironment{proposition}{\refstepcounter{proposition} \emph{Proposition \theproposition.}}{}
\newenvironment{proposition*}[1]{\emph{Proposition #1.}}{}
\newcounter{theorem}
\newenvironment{theorem}{\refstepcounter{theorem} \emph{Theorem \thetheorem.}}{}
\newenvironment{theorem*}[1]{\emph{Theorem #1.}}{}
\renewenvironment{proof}{\emph{Proof.}}{$\square$}
\newenvironment{proof*}[1]{\emph{#1.}}{$\square$}

\DeclareMathOperator{\First}{First}
\DeclareMathOperator{\Last}{Last}

\aicAUTHORdetails{%
  title = {Leaper Tours},
  author = {Nikolai Beluhov},
  plaintextauthor = {Nikolai Beluhov},
  keywords = {Leaper graphs, Hamiltonian cycles},
}

\aicEDITORdetails{%
   year={2022},
   number={4},
   received={28 April 2021},
   revised={28 December 2021},
   published={1 June 2022},
   doi={10.19086/aic.2022.4},
}

\makeatletter
\renewcommand{\toc@lastpage}{32}
\makeatother

\begin{document}

\begin{frontmatter}[classification=text]

\author[nb]{Nikolai Beluhov}

\begin{abstract}
Let $p$ and $q$ be positive integers. The $(p, q)$-leaper $L$ is a generalised knight which leaps $p$ units away along one coordinate axis and $q$ units away along the other. Consider a free $L$, meaning that $p + q$ is odd and $p$ and $q$ are relatively prime. We prove that $L$ tours the board of size $4pq \times n$ for all sufficiently large positive integers $n$. Combining this with the recently established conjecture of Willcocks which states that $L$ tours the square board of side $2(p + q)$, we conclude that furthermore $L$ tours all boards both of whose sides are even and sufficiently large. This, in particular, completely resolves the question of the Hamiltonicity of leaper graphs on sufficiently large square boards.
\end{abstract}

\end{frontmatter}

\section{Introduction} \label{intro}

The \emph{$(p, q)$-leaper} $L$ is a fairy chess piece that generalises the knight. On each move, it leaps $p$ units away along one coordinate axis and $q$ units away along the other. For the knight, $p = 1$ and $q = 2$. Other named leapers include the wazir ($p = 0$ and $q = 1$), giraffe ($p = 1$ and $q = 4$), zebra ($p = 2$ and $q = 3$), flamingo ($p = 1$ and $q = 6$), and antelope ($p = 3$ and $q = 4$).

Let $A$ be a board of size $m \times n$. We give the number of rows first and the number of columns second, so that $A$ has $m$ rows, $n$ columns, and $mn$ cells.

A \emph{tour} of $L$ on $A$ is a sequence of moves of $L$ on $A$ such that $L$ visits every cell of $A$ exactly once and with its final move returns where it started. This is also known as a \emph{closed tour}. The definition of an \emph{open tour} is the same, except that $L$ starts and finishes on different cells.

(We use ``tour'' and ``closed tour'' synonymously as we touch upon open tours only rarely. We do not consider an empty sequence of moves on the board of size $1 \times 1$ or a two-move back-and-forth by the wazir on a domino board to be tours. Thus a tour always visits at least three cells over the course of at least three moves.)

When $p + q$ is even, $L$ cannot tour any board since it only visits cells of the same colour in the traditional chessboard colouring. Similarly, when $d = \gcd(p, q) \ge 2$, once again $L$ cannot tour any board because its moves never change the coordinates of its cell modulo $d$.

The leapers $L$ such that $p + q$ is odd and $p$ and $q$ are relatively prime are known as \emph{free}. Thus free leapers are the only ones which could potentially have tours. Suppose, from now on, that $L$ is free.

When $A$ is too narrow, $L$ might not be able to tour it because it does not have enough room for manoeuvres. So for a tour we need both of $m$ and $n$ to be reasonably large. However, even then we could still run into problems: Since each move of $L$ changes the colour of its cell, a tour requires an even number of moves. Thus $A$ must also have an even number of cells.

It appears plausible that these two are the only obstacles which could possibly prevent $L$ from touring a board. Or, in other words, that $L$ tours all sufficiently large boards with even area. Here, by ``sufficiently large'' we mean all boards both of whose sides are greater than or equal to some fixed positive integer depending on $L$.

Soon we will discuss this conjecture in more detail. First, however, let us take care of some technicalities.

We can easily show that the wazir tours $A$ if and only if $m \ge 2$, $n \ge 2$, and $mn$ is even. Thus for the wazir the mere existence of tours is not a very difficult question. The behaviour of the wazir differs from the behaviour of all other free leapers in other important ways as well. 

A \emph{skew} leaper is one for which $p \neq 0$, $q \neq 0$, and $p \neq q$. Thus we can exclude the wazir by specifying that we will consider only skew free leapers. Suppose, from now on, that $L$ is indeed both skew and free.

Finally, since swapping $p$ and $q$ has no effect on $L$, let us suppose for concreteness that $p < q$.

Or, in summary: From now on throughout the rest of the present work, to us $p$ and $q$ will be positive integers with $p < q$ such that $p + q$ is odd and $p$ and $q$ are relatively prime, and $L$ will be the skew free $(p, q)$-leaper.

Some of the most natural questions we can ask about leaper tours are as follows.

\medskip

\begin{question*}{\textbf{A}} Does $L$ tour some board? \end{question*}

\medskip

\begin{question*}{\textbf{B}} Does $L$ tour arbitrarily large boards? \end{question*}

\medskip

\begin{question*}{\textbf{C}} Does $L$ tour all sufficiently large boards with even area? \end{question*}

\medskip

\begin{question*}{\textbf{D}} What are all boards which $L$ tours? \end{question*}

\medskip

We can consider each one of these questions for a concrete leaper, or for some infinite family of leapers, or for all skew free leapers.

The most closely studied concrete leaper is, of course, the knight.

The two infinite families of skew free leapers which have received the most attention are, by a wide margin, $p = 1$ (beginning with the knight, giraffe, and flamingo) and $q = p + 1$ (beginning with the knight, zebra, and antelope). We call them the two \emph{exceptional} families. Experience shows that for many problems involving leapers it makes sense to consider them first as they are likely to contain the lowest-hanging fruit.

Another natural restriction in Questions \textbf{A}--\textbf{D} would be to consider square boards only, instead of arbitrary rectangular boards.

Often a solution to Question \textbf{A} will also settle Question \textbf{B} with little additional effort. Indeed, suppose that we have managed somehow to find a tour $T$ of $L$ on $A$. Take any board $B$ of the form $am \times bn$ and tile it with copies of $A$ each one of which also contains a copy of $T$. Then we might be able to construct a tour of $B$ by deleting some small number of moves and inserting new ones in their place so as to stitch all such copies of $T$ together into one single tour. That would resolve Question \textbf{B} as well.

We cannot guarantee beforehand that such an approach would work. Ordinarily, though, it does. Thus we will talk more about Question \textbf{A} and less about Question \textbf{B}, keeping in mind this connection between them.

Using our hypothetical tour $T$ of $A$ as a building block, we can at best hope to construct new tours on those boards $B$ whose area is a multiple of $mn$. So for Question \textbf{C} we need either multiple building blocks which furthermore can be stitched together in different combinations, or else a completely novel strategy for tour construction. Thus we should expect Question \textbf{C} to be distinctly more difficult than Questions \textbf{A} and \textbf{B}. This has indeed been borne out in practice.

The answer to Question \textbf{D} probably will not be tractable in the general case. Thus Question \textbf{C} seems like a reasonable compromise between what we would like to know and what we are likely to be able to learn.

We go on to a brief historical overview. We pay close attention to the generality of the results: Do they apply to a single leaper, or to some infinite family, or to all skew free leapers? Furthermore, what boards do they apply to?

(Though note that all instances of the words ``first'', ``earliest'', and ``oldest'' in what follows should be understood to include the caveat ``so far as the author has been able to determine at the time of writing''.)

The oldest surviving knight tour, on the $8 \times 8$ board, is found in a ninth century manuscript of Sha\d{t}ranj champion al-`Adl\={\i} ar-R\=um\=\i. (Though some open tours of the $4 \times 8$ board have been recovered from slightly earlier works as well.)

Since then, knight tours have given rise to a vast amount of literature. We point readers to Jelliss's \cite{J2} for an in-depth discussion of the subject as well as a comprehensive bibliography.

The earliest mathematically sophisticated proof of nonexistence for the knight was given in 1877 by Flye Sainte-Marie: The knight does not tour any board of the form $4 \times n$. \cite{FSM} (Strictly speaking, Flye Sainte-Marie considered only $n = 8$. However, his argument is fully general. By contrast, open tours do exist when $n = 3$ or $n \ge 5$.)

A positive answer to Question \textbf{C} for the knight was only spelled out explicitly and rigorously in 1991 by Schwenk. \cite{S} In fact, the same work resolves the full Question \textbf{D} for the knight as well. Schwenk himself points out in his introduction how surprising it is that this was not done much earlier.

The first tours of skew free leapers other than the knight are due to Frost and appeared in the 1886 book \cite{F}. Frost constructed one tour each of the giraffe and the zebra on the board of size $10 \times 10$. Many more tours of concrete leapers on concrete boards followed in the twentieth century. We refer readers to volume 10 of \cite{J2} for a detailed account with reproductions of the tours.

The earliest mathematically sophisticated proof of nonexistence for a concrete skew free leaper other than the knight was given in 1976 by Jelliss: The giraffe does not tour the board of size $8 \times 8$. \cite{JW} As Jelliss himself remarks in \cite{J1}, his argument shows also that the giraffe does not tour any board of the form $8 \times n$. (Though open tours do exist on some of these boards.)

Dawson showed in 1928 that, when $L$ belongs to the exceptional family $p = 1$, there exists an open tour of $L$ from corner to corner on the board of size $(q + 1) \times 2q$. \cite{Do} This was the first construction of either open or closed tours for an infinite family of skew free leapers.

Willcocks conjectured in 1976 that the smallest square board which $L$ tours is the one of side $2(p + q)$. \cite{JW} (This is also true for the wazir.) There are two distinct parts to Willcocks's conjecture: That tours do not exist on smaller square boards and that $L$ does tour the square board of side $2(p + q)$. As we are about to see, the two parts differ greatly in difficulty. Frost's tours confirm the existence part of the conjecture for the giraffe and the zebra, and Willcocks himself constructed tours confirming it for the $(2, 5)$-leaper and the antelope.

The special case of the exceptional family $p = 1$ was settled by Dejter in 1988. \cite{De} In fact, the same work resolves Question \textbf{B} for these leapers as well. This was the earliest construction of closed tours for an infinite family of leapers.

One more important milestone is Jelliss's \cite{J1} which collects and systematises more or less all knowledge about leaper tours obtained by the time when it was written in 1993.

Then, in 1994, Knuth published the seminal work \cite{Kn} containing a wealth of original results. We cover only some of the major highlights. Knuth gave a full proof of the nonexistence part of Willcocks's conjecture together with constructions verifying its existence part for all skew free leapers in the two exceptional families $p = 1$ and $q = p + 1$. (The former independently from Dejter.) He also generalised Flye Sainte-Marie's and Jelliss's proofs to show that a skew free leaper $L$ can never tour a board of height~$2q$.

Thus in particular Knuth authored the first mathematically sophisticated proofs of nonexistence for all skew free leapers and found novel constructions of closed tours regarding infinite families of leapers.

Ong established a number of further impossibility results in 2001. \cite{O} One notable among them is that, in the special case when $L$ belongs to the exceptional family $q = p + 1$, it does not tour any board of height $2(ap + b)$ with $1 \le a \le b \le p$. One more is that an arbitrary skew free leaper $L$ can never tour a board of height $p + q + 2a - 1$ with $1 \le a \le p$. (The arguments are, as before, specifically about closed tours and do not address open ones.)

Kam\v{c}ev, in 2014, solved the special case of Question \textbf{C} for skew free leapers $L$ in the exceptional family $p = 1$ on square boards. \cite{Ka}

We said in the beginning that the answer to Question \textbf{C} is most probably positive for all skew free leapers. This conjecture was stated explicitly for the first time in 2012 by Erde. \cite{E}

\medskip

\begin{conjecture*}{\textbf{E}} (Even area.) For all sufficiently large boards $A$, the skew free leaper $L$ tours $A$ if and only if the area of $A$ is even. \end{conjecture*}

\medskip

In 2017, the present author proved the existence part of Willcocks's conjecture for all skew free leapers. \cite{B} The proof in particular supplied the first positive answer to Question \textbf{A} in full generality. The same work handles Question \textbf{B} as well, with the help of the strategy we outlined in our initial discussion of Questions \textbf{A}--\textbf{D}.

This sums up the current state of the art: Questions \textbf{A} and \textbf{B} have been resolved for all free leapers. Question \textbf{C}, on the other hand, has been settled completely only for the wazir and the knight. (Though note that, when $p$ and $q$ are reasonably small, we can likely tackle Question \textbf{C} for $L$ without too much difficulty by constructing a suitable set of building blocks with the help of modern computational technology.) Furthermore, Question \textbf{C} has also been settled in the special case when $L$ belongs to the exceptional family $p = 1$ with a restriction to square boards.

Our main goal in the present work will be to prove the following theorems.

\medskip

\begin{theorem} \label{4pq} Let $L$ be a skew free $(p, q)$-leaper. Then $L$ tours all boards of width $4pq$ and sufficiently large height. \end{theorem}

\medskip

Our proof of Theorem \ref{4pq} is completely self-contained. Thus it answers Question \textbf{A} in the affirmative in a new way independent from the one in \cite{B}.

Then we use the tours of Willcocks's conjecture and the tours of Theorem \ref{4pq} together as building blocks so as to construct tours of $L$ on many other boards, too.

\medskip

\begin{theorem} \label{main} Let $L$ be a skew free leaper. Then $L$ tours all boards both of whose sides are even and sufficiently large. \end{theorem}

\medskip

Theorem \ref{main} is a significant advance towards the full even area conjecture.

One immediate corollary of Theorem \ref{main} is worth mentioning as well.

\medskip

\begin{theorem} \label{square} Let $L$ be a skew free leaper. Then for all sufficiently large positive integers $n$ we have that $L$ tours the square board of side $n$ if and only if $n$ is even. \end{theorem}

\medskip

Theorem \ref{square} completely resolves Question \textbf{C} on square boards for all skew free leapers. Thus, in particular, it also confirms this special case of the even area conjecture.

\section{Preliminaries} \label{prelim}

We continue with the formalisation of cells, boards, and leaper tours.

We model a cell as an ordered pair of integers and a board as the Cartesian product of two integer intervals. Thus, in particular, we model a board as a set of cells. By convention, our coordinates increase from left to right and upwards. Note also that we do not require the rows and columns of a board to be numbered beginning with some predetermined starting point such as zero or one.

We define the parity of cell $(x, y)$ to be the parity of the sum $x + y$. Thus even and odd cells form the exact same pattern as the black and white cells in the traditional chessboard colouring.

We talk about leaper tours in the language of graph theory. In terms of coordinates, two cells $(x', y')$ and $(x'', y'')$ are the endpoints of an edge of $L$ when $\{|x' - x''|, |y' - y''|\} = \{p, q\}$. Then we define the \emph{leaper graph} $G$ of $L$ on a board $A$ as follows: The vertices of $G$ are the cells of $A$ and the edges of $G$ are all edges of $L$ on $A$. Note that $G$ is unoriented even though moves in chess have definite directions.

Thus a tour of $L$ on $A$ corresponds to a Hamiltonian cycle in $G$. Similarly, an open tour becomes a Hamiltonian path in $G$. From now on, we will use the terms ``tour'' and ``Hamiltonian cycle'' synonymously. We do not consider open tours and Hamiltonian paths outside of the introduction, and so when we say that a graph is Hamiltonian we will always mean that it admits a closed tour.

A \emph{pseudotour} of $L$ on $A$ is a collection of cycles of $L$ on $A$ visiting each cell of $A$ exactly once. Equivalently, a pseudotour is a subgraph of the leaper graph of $L$ on $A$ in which every cell of $A$ has degree two. Thus locally a pseudotour looks just like a tour. Globally, however, a pseudotour need not be connected. (A pseudotour can also be described as a \emph{two-factor} or a \emph{cycle cover}.)

We conclude with a couple of remarks on notation.

We use both unary and binary $\bmod$. Thus $a \equiv b \pmod d$ means that $d$ divides $a - b$, whereas $a \bmod d$ stands for the unique nonnegative integer $r$ with $r < d$ such that $d$ divides $a - r$.

We also denote the edge joining vertices $u$ and $v$ by either $uv$ or $u$---$v$, whichever one reads more clearly in the specific setting at hand.

\section{Tours in Projection Graphs I} \label{proj-i}

Let $I$ be an integer interval. The \emph{projection graph} $\Pi(p, q, I)$ has vertex set $I$ and edges $uv$ for all $u$ and $v$ in $I$ that differ by either $p$ or $q$. We say that edge $uv$ is \emph{short} when $|u - v| = p$ and \emph{long} otherwise, when $|u - v| = q$. When only the size $n = |I|$ of $I$ matters, we write simply $\Pi(p, q, n)$. In such cases, for concreteness we also assume that $I = [0; n - 1]$ unless otherwise stated.

\begin{figure}[ht] \centering \includegraphics{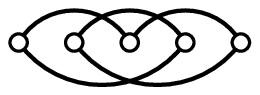} \caption{} \label{f:proj} \end{figure}

For example, Figure \ref{f:proj} shows $\Pi(2, 3, 5)$.

We call these objects projection graphs because they are what we get when we project leaper graphs onto the coordinate axes. They will help us construct and study leaper pseudotours. Before that, however, we must learn to construct some tours of the projection graphs themselves.

Ginzboorg and Niemi in \cite{GN} study the closely related family of graphs obtained when we orient all short edges one way and all long edges the opposite way. They also point out the connection with leaper graphs and establish a result equivalent to our Lemma \ref{pst} below.

We introduce some notations so that we can describe our constructions clearly and succinctly.

Let $I = [a; b] = \{u \mid a \le u \le b\}$. We write $\First(k, I)$ for the integer interval formed by the first $k$ elements of $I$, namely $[a; a + k - 1]$. Similarly, we write $\Last(k, I)$ for the integer interval formed by the last $k$ elements of $I$, namely $[b - k + 1; b]$. When $k > |I|$, we set $\First(k, I) = \Last(k, I) = I$.

When we say something like ``partition $I$ into $k$ subintervals $I_1$, $I_2$, \ldots, $I_k$ of sizes $s_1$, $s_2$, \ldots, $s_k$, respectively'', we always mean also that $I_1$, $I_2$, \ldots, $I_k$ lie on the number line in this order from left to right. Thus explicitly $I_1 = \First(s_1, I)$, $I_2 = [a + s_1; a + s_1 + s_2 - 1]$, $I_3 = [a + s_1 + s_2; a + s_1 + s_2 + s_3 - 1]$, \ldots, $I_k = \Last(s_k, I)$.

We say that $a + i$ is the element at \emph{position} $i$ in $I$. Consider two integer intervals $I'$ and $I''$ of the same size. Then we say that $u' \in I'$ and $u'' \in I''$ are \emph{corresponding} elements of $I'$ and $I''$ when they occupy identical positions in $I'$ and $I''$, respectively.

We write $I + k$ for the integer interval $[a + k; b + k]$. Let $I' + k = I''$ with $k \in \{\pm p, \pm q\}$. Then we write $I' \circeq I''$ for the set of all edges of $\Pi(p, q, \mathbb{Z})$ which join corresponding elements of $I'$ and $I''$. There are $|I'| = |I''|$ such edges, all of which are short when $|k| = p$ and all of which are long when $|k| = q$. We call a set of edges of this form a \emph{pencil}.

\medskip

\begin{lemma} \label{pst} The projection graph $\Pi(p, q, p + q)$ is Hamiltonian. \end{lemma} 

\medskip

\begin{proof} This is something of an understatement. In fact, $\Pi(p, q, p + q)$ consists of a single cycle.

Perhaps the simplest argument is as follows. Let $u_i = pi \bmod (p + q)$. Then every vertex of $\Pi(p, q, p + q)$ has the form $u_i$ for some $i$ because $p$ and $q$ are relatively prime. On the other hand, since addition of $p$ and subtraction of $q$ coincide modulo $p + q$, we obtain furthermore that the neighbours of $u_i$ in $\Pi(p, q, p + q)$ are $u_{i - 1}$ and $u_{i + 1}$ for all $i$. \end{proof}

\medskip

We can construct larger tours using $\Pi(p, q, p + q)$ as a building block.

\medskip

\begin{lemma} \label{pmt} The projection graph $\Pi(p, q, k(p + q))$ is Hamiltonian for all $k$. \end{lemma} 

\medskip

\begin{proof} Partition an integer interval of size $k(p + q)$ into subintervals $I_1$, $I_2$, \ldots, $I_k$, each of size $p + q$, and let $C_i = \Pi(p, q, I_i)$. By the proof of Lemma \ref{pst}, each $C_i$ is a cycle.

To stitch all such cycles together into a tour, let $u_i$, $v_i$, and $w_i$ be the vertices of $I_i$ at positions $\lfloor (q - p)/2 \rfloor$, $\lfloor (p + q)/2 \rfloor$, and $\lfloor (3p + q)/2 \rfloor$, respectively. Then for all $i$ with $1 \le i < k$ delete the pair of short edges $v_iw_i$ and $u_{i + 1}v_{i + 1}$ and replace them with the pair of long edges $v_iu_{i + 1}$ and $w_iv_{i + 1}$. \end{proof}

\medskip

We proceed to develop one more technique for the construction of larger tours out of smaller ones.

Let $I' = [a'; b']$ and $I'' = [a''; b'']$ be two disjoint integer intervals of the same size $s$. Suppose, for concreteness, that $I'$ and $I''$ lie in this order on the number line, so that $I'' = I' + t$ for some positive integer $t$ with $t \ge s$.

Let $J = [a'; b'']$ be the smallest integer interval which includes both of $I'$ and $I''$. An \emph{extension} between $I'$ and $I''$ is a subgraph $\mathcal{E}$ of $\Pi(p, q, J)$ such that: (a) It is the disjoint union of $s$ paths each one of which connects two corresponding vertices of $I'$ and $I''$; and (b) The paths of $\mathcal{E}$ visit all vertices of $J$. Thus the vertices of $I'$ and $I''$ have degree one in $\mathcal{E}$, and all other vertices of $J$ have degree two in $\mathcal{E}$. We call $s$ the \emph{width} of $\mathcal{E}$ and we call $t$ its \emph{length}.

We can concatenate together several extensions of the same width, as follows. Consider an integer interval of size $s + kt$ partitioned into $2k + 1$ subintervals $I_1$, $J_1$, $I_2$, $J_2$, \ldots, $I_{k + 1}$ of sizes $s$, $t - s$, $s$, $t - s$, \ldots, $s$, respectively. Then for all $i$ with $1 \le i \le k$ construct a translation copy of $\mathcal{E}$ on vertex set $I_i \cup J_i \cup I_{i + 1}$. We call the union of all such copies a \emph{chain} of $\mathcal{E}$. Clearly, a chain of $k$ instances of $\mathcal{E}$ is itself an extension of the same width $s$ and length $kt$.

We will return to extensions in a moment.

Let $T$ be a tour of $\Pi(p, q, I)$. Consider any partitioning of the edges of $T$ into two subsets $T'$ and $T''$. For each $i = 0$, $1$, and $2$, let $V_i$ be the set of all vertices of $\Pi(p, q, I)$ with degree $2 - i$ in $T'$ and degree $i$ in~$T''$.

We say that $T'$ and $T''$ form a \emph{split} of $T$ when: (a) All three of $V_0$, $V_1$, and $V_2$ are integer intervals; and (b) These integer intervals lie on the number line in this order from left to right. We call the size of $V_1$ the \emph{width} of the split. We allow $V_0$ and $V_2$ to be empty, but not $V_1$. Thus the width of a split is always positive.

The point of these definitions is as follows.

\medskip

\begin{lemma} \label{es} Suppose that there exists an extension of width $s$ and length $t$. Suppose, furthermore, that $\Pi(p, q, n)$ is Hamiltonian and one of its tours admits a split of width $s$. Then $\Pi(p, q, n + kt)$ is Hamiltonian as well for all $k$. \end{lemma} 

\medskip

\begin{proof} Let $\mathcal{E}$ be that extension, let $T$ be that tour, and let $T'$ and $T''$ be the two parts of that split of $T$. Then we can build the required tour of $\Pi(p, q, n + kt)$ out of two translation copies of $T'$ and $T''$ connected by a chain of $k$ instances of $\mathcal{E}$. \end{proof}

\medskip

\begin{figure}[ht!] \centering \begin{tabular} {>{\centering}m{0.95\textwidth}} \includegraphics{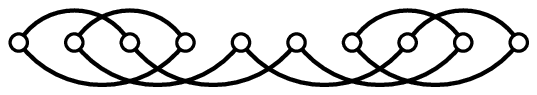} \tabularnewline \bigbreak \tabularnewline \includegraphics{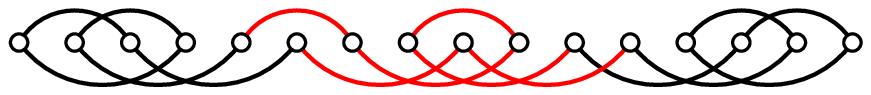} \end{tabular} \caption{} \label{f:es} \end{figure}

Figure \ref{f:es} shows an example with $s = 2$, $t = 6$, and $k = 1$. The original tour is shown above and the extended tour is shown below, with the extension's edges highlighted.

We already know by Lemma \ref{pmt} that $\Pi(p, q, n)$ is Hamiltonian for infinitely many $n$. Equipped with the technique of extensions and splits, we proceed to show that in fact it is Hamiltonian for all sufficiently large $n$.

First let us deal with the case of the knight.

\medskip

\begin{lemma} \label{pkt} The projection graph $\Pi(1, 2, n)$ is Hamiltonian if and only if $n \ge 3$. Furthermore, for all such $n$ it admits a unique tour. That tour consists of the short edges $0$---$1$ and $(n - 2)$---$(n - 1)$ together with all long edges in the graph. \end{lemma} 

\medskip

\begin{proof} Consider all vertices of $\Pi(1, 2, n)$ one by one, working from left to right, and all the while keeping careful track of which edges are forced to belong to the tour and which ones are forced to be outside of it. \end{proof}

\medskip

For all skew free leapers other than the knight, $2p \neq q$. The cases when $2p < q$ and $2p > q$ behave very differently and we consider them separately. The former case is much simpler and so we take care of it first.

\medskip

\begin{lemma} \label{ee} Suppose that $2p < q$. Then there exists an extension of width $2p$ and length $2q$. \end{lemma} 

\medskip

\begin{proof} For concreteness, let $I' = [0; 2p - 1]$ and $I'' = I' + 2q = [2q; 2(p + q) - 1]$. Let also $I_\text{Mid} = I' + q = I'' - q$.

To begin with, consider the union of all paths of the form $(u - q)$---$u$---$(u + q)$ with $u \in I_\text{Mid}$. These paths are pairwise disjoint, they connect all pairs of corresponding vertices of $I'$ and $I''$, and they visit all vertices of $I_\text{Mid}$. However, they do not quite form an extension because they leave all vertices of the subintervals $J' = [2p; q - 1]$ and $J'' = [2p + q; 2q - 1]$ unvisited.

\begin{figure}[ht] \centering \includegraphics{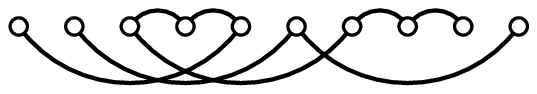} \caption{} \label{f:ee} \end{figure} 

We patch things up as follows. For each vertex $u$ in $\First(p, J')$, find the smallest positive integer $i$ with $v = u + pi \in I_\text{Mid}$. Then delete edge $v$---$(v + q)$ and replace it with the path $v$---$(v - p)$---$(v - 2p)$---$\cdots$---$u$---$(u + q)$---$(u + p + q)$---$(u + 2p + q)$---$\cdots$---$(v + q)$. The net result of all such rewirings over all vertices $u$ in $\First(p, J')$ will be an extension between $I'$ and $I''$. \end{proof}

\medskip

For example, Figure \ref{f:ee} shows the resulting extension for the giraffe.

\medskip

\begin{lemma} \label{se} Suppose that $2p < q$. Then every tour in a projection graph admits a split of width $2p$. \end{lemma} 

\medskip

\begin{proof} Let $T$ be a tour of $\Pi(p, q, I)$. Since all vertices of $\First(p, I)$ have degree two in $\Pi(p, q, I)$, all short edges of the form $u$---$(u + p)$ with $u \in \First(p, I)$ must necessarily belong to $T$. Let $T'$ be the set of all such edges and let $T''$ be the set of all other edges of $T$. Then $T'$ and $T''$ satisfy the definition of a split with $V_1 = \First(2p, I)$. \end{proof}

\medskip

\begin{lemma} \label{pte} Suppose that $2p < q$. Then $\Pi(p, q, n)$ is Hamiltonian for all sufficiently large $n$. \end{lemma} 

\medskip

\begin{proof} By Lemmas \ref{pmt}, \ref{es}, \ref{ee}, and \ref{se} we obtain that $\Pi(p, q, n)$ is Hamiltonian for every $n$ of the form $n = (p + q)i + 2qj$ where $i$ is a positive integer and $j$ is a nonnegative integer. \end{proof}

\medskip

This settles the case when $2p < q$. The case when $2p > q$ is considerably more difficult, and to it we devote the entire next section.

\section{Tours in Projection Graphs II} \label{proj-ii}

Suppose that $2p > q$. Let $d = q - p$, $r = p \bmod d = q \bmod d$, and $\alpha = \lfloor q/d \rfloor$. Thus $\alpha \ge 2$, $p = (\alpha - 1)d + r$, and $q = \alpha d + r$.

First we construct extensions for this case.

\medskip

\begin{lemma} \label{ed} Suppose that $\max\{2d, q - 2d\} \le s < q$. Then there exists an extension of width $s$ and length $2q$. \end{lemma} 

\medskip

\begin{proof} Let $I' = [0; s - 1]$, $I'' = I' + 2q = [2q; s + 2q - 1]$, and $I_\text{Mid} = I' + q = I'' - q$. As in the proof of Lemma \ref{ee}, consider initially the paths $(u - q)$---$u$---$(u + q)$ for all $u$ in $I_\text{Mid}$. The sets of vertices that they leave unvisited are $J' = [s; q - 1]$ and $J'' = [s + q; 2q - 1]$.

We patch things up as follows. For each vertex $u$ in $J'$, set $v = u + p$ when $u \in \First(d, J')$ and $v = u - p$ otherwise, when $u \not \in \First(d, J')$, so that $v \in I_\text{Mid}$ in the former case and $v \in I'$ in the latter one. Then delete edge $v$---$(v + q)$ and replace it with the path $v$---$u$---$(u + q)$---$(v + q)$. The net result of all such rewirings over all vertices $u$ in $J'$ will be an extension between $I'$ and $I''$. \end{proof}

\medskip

\begin{figure}[ht] \centering \includegraphics[scale=0.85]{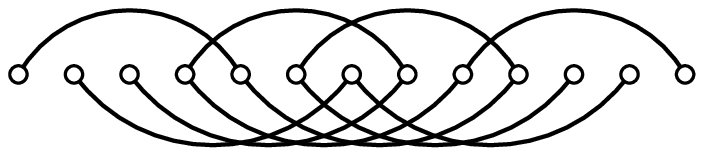} \caption{} \label{f:ed} \end{figure}

For example, Figure \ref{f:ed} shows the resulting extension with $s = 3$ for the $(4, 5)$-leaper.

The proof of Lemma \ref{ed} goes through also with the weaker constraints $\max\{d, q - 2d\} \le s \le q$. The reason we use stronger bounds will become clear in the proof of Lemma \ref{swb}.

We go on to construct some tours which admit appropriate splits.

Let $k$ be a positive integer such that $2k \le \alpha$ and let $I$ be an integer interval of size $2kq$. Partition $\First(2kp, I)$ into $2k$ subintervals $P_1$, $P_2$, \ldots, $P_{2k}$ of size $p$ each, and also partition $I$ itself into $2k$ subintervals $Q_1$, $Q_2$, \ldots, $Q_{2k}$ of size $q$ each.

We define $\mathcal{B}_k$ to be the subgraph of $\Pi(p, q, I)$ whose edge set is the union of the pencils $P_{2i - 1} \circeq P_{2i}$ and $Q_{2i - 1} \circeq Q_{2i}$ for all $i$ with $1 \le i \le k$.

Clearly, all vertices of $\mathcal{B}_k$ in $\Last(2kd, I)$ have degree one and all other vertices of $\mathcal{B}_k$, namely the ones in $\First(2kp, I)$, have degree two. Therefore, $\mathcal{B}_k$ is the disjoint union of $kd$ paths and, possibly, some cycles.

We proceed to examine the structure of $\mathcal{B}_k$ into somewhat deeper detail. To that end, partition $\Last(2kd, I)$ into $2k$ subintervals $D_1$, $D_2$, \ldots, $D_{2k}$ of size $d$ each.

\medskip

\begin{lemma} \label{bookend} A path of $\mathcal{B}_k$ connects the vertex at position $i$ in $D_j$ to the vertex at position $(i + r) \bmod d$ in $D_{(2k + 1) - j}$ for all $i$ and $j$ with $0 \le i < d$ and $1 \le j \le k$. Furthermore, $\mathcal{B}_k$ does not contain any cycles. \end{lemma}

\medskip

\begin{proof} We begin with $j = 1$. Let $u$ be the vertex at position $i$ in $D_1$. Initially, the path starting from $u$ alternates between long and short edges. The first two steps along it bring us to $u - d$, then the next two bring us to $u - 2d$, and so on. Eventually, after some number of steps $2\ell$, we will reach a vertex $u - d\ell$ in $\Last(d, Q_{2k - 1})$. One more long edge brings us to $u - d\ell + q$. This is the vertex at position $(i + r) \bmod d$ in $D_{2k}$, and there the journey ends.

When $2 \le j \le k$, we proceed by induction on $k$. Our base case is $k = 1$, when there is nothing to prove. For the induction step, suppose that $k \ge 2$. Observe that $\mathcal{B}_k$ is the disjoint union of one copy of $\mathcal{B}_{k - 1}$ on vertices $\First(2(k - 1)q, I)$ together with the pencils $P_{2k - 1} \circeq P_{2k}$ and $Q_{2k - 1} \circeq Q_{2k}$. The vertices of that copy of $\mathcal{B}_{k - 1}$ with degree one are $\Last(2(k - 1)d, Q_{2k - 2})$. Partition them into subintervals $D'_1$, $D'_2$, \ldots, $D'_{2k - 2}$ of size $d$ each.

Let $u$ be the vertex at position $i$ in $D_j$, with $2 \le j \le k$. Then the path starting from $u$ begins with one long edge and one short edge which connect $u$ to the corresponding vertex $u'$ in subinterval $D'_{j - 1}$. By our induction hypothesis, a path of $\mathcal{B}_{k - 1}$ then connects $u'$ to the vertex $u''$ at position $(i + r) \bmod d$ in subinterval $D'_{2k - j}$. Finally, one more short edge followed by one more long edge connect $u''$ to the corresponding vertex in subinterval $D_{(2k + 1) - j}$, as needed.

Since the above accounts for all edges of $\mathcal{B}_k$, additionally we see that it does not contain any cycles.~\end{proof}

\medskip

We continue with the construction of the tours themselves. The cases when $\alpha$ is even versus odd behave very differently and we consider them separately. The case when $\alpha$ is even is much simpler and so we take care of it first.

\medskip

\begin{lemma} \label{sde} Suppose that $2p > q$ and $\alpha$ is even. Then $\Pi(p, q, \alpha(p + q) + d)$ is Hamiltonian. \end{lemma} 

\medskip

\begin{proof} Let $I$ be an integer interval of size $\alpha(p + q) + d$. Let $k = \alpha/2$ and partition $I$ into subintervals $J'$, $D_1$, $D_2$, \ldots, $D_{2k + 1}$, $J''$ of sizes $2kp$, $d$, $d$, \ldots, $d$, $2kp$, respectively.

Place one translation copy $\mathcal{B}'$ of $\mathcal{B}_k$ so that its vertices of degree one are $D_1 \cup D_2 \cup \cdots \cup D_{2k}$ and its vertices of degree two are $J'$; as well as one reflection copy $\mathcal{B}''$ of $\mathcal{B}_k$ so that its vertices of degree one are $D_2 \cup D_3 \cup \cdots \cup D_{2k + 1}$ and its vertices of degree two are $J''$.

By this point in our construction, all vertices of $I$ have degree two except for the ones in $D_1$ and $D_{2k + 1}$. To fix that, we add in also the two pencils $\First(d - r, D_1) \circeq \Last(d - r, D_{2k + 1})$ and $\Last(r, D_1) \circeq \First(r, D_{2k + 1})$. Observe that these pencils join each vertex at position $i$ in $D_1$ to the vertex at position $(i + r) \bmod d$ in $D_{2k + 1}$.

We are left to verify that our construction does indeed yield a tour, as opposed to a disjoint union of several cycles.

Consider the vertex of $D_1$ at position $i$. By Lemma \ref{bookend}, a path in $\mathcal{B}'$ leads from it to the vertex at position $(i + r) \bmod d$ in $D_{2k}$. Then, again by Lemma \ref{bookend}, another path in $\mathcal{B}''$ leads from there to the vertex of $D_3$ at position $i$ when $k \ge 2$ and $i + 2r$ when $k = 1$. We continue on in this way, travelling alternatively along paths in $\mathcal{B}'$ and $\mathcal{B}''$, and visiting successively $D_{2k - 2}$, $D_5$, \ldots, $D_2$, $D_{2k + 1}$. One final edge brings us from $D_{2k + 1}$ to the vertex at position $(i + (-1)^{k + 1}r) \bmod d$ in $D_1$.

Observe that $r$ and $d$ are relatively prime. (This also holds in the case of $q = p + 1$, when $r = 0$ and $d = 1$.) Thus over the course of $d$ iterations of the above itinerary we will visit all vertices of $D_1$, returning where we started, and meanwhile traversing the entirety of $\mathcal{B}'$ and $\mathcal{B}''$ together with all of the newly added edges between $D_1$ and $D_{2k + 1}$. Therefore, our construction does indeed yield a tour, as needed. \end{proof}

\medskip

\begin{figure}[ht] \centering \includegraphics[scale=0.45]{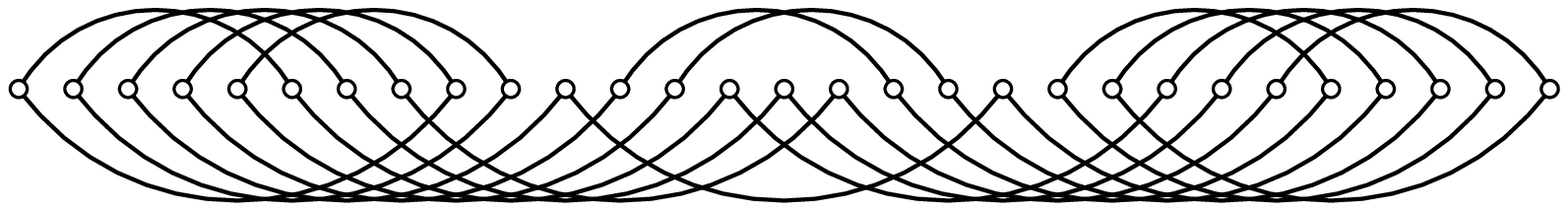} \caption{} \label{f:sde} \end{figure}

Figure \ref{f:sde} shows an example with $p = 5$ and $q = 8$, when $d = 3$, $r = 2$, and $\alpha = 2$.

We go on to the case when $\alpha$ is odd.

\medskip

\begin{lemma} \label{sdo} Suppose that $2p > q$ and $\alpha$ is odd. Then $\Pi(p, q, \alpha(p + q) + 2d)$ is Hamiltonian. \end{lemma} 

\medskip

\begin{proof} Let $I$ be an integer interval of size $\alpha(p + q) + 2d$. Let $k = \lfloor \alpha/2 \rfloor$ and partition $I$ into subintervals $J'$, $D'_1$, $D'_2$, \ldots, $D'_{2k}$, $J_\text{Mid}$, $D''_1$, $D''_2$, \ldots, $D''_{2k}$, $J''$ of sizes $2kp$, $d$, $d$, \ldots, $d$, $(2k + 3)d + 2r$, $d$, $d$, \ldots, $d$, $2kp$, respectively.

As in the proof of Lemma \ref{sde}, place one translation copy $\mathcal{B}'$ of $\mathcal{B}_k$ so that its vertices of degree one are $D'_1 \cup D'_2 \cup \cdots \cup D'_{2k}$ and its vertices of degree two are $J'$; as well as one reflection copy $\mathcal{B}''$ of $\mathcal{B}_k$ so that its vertices of degree one are $D''_1 \cup D''_2 \cup \cdots \cup D''_{2k}$ and its vertices of degree two are $J''$.

Then furthermore partition $J_\text{Mid}$ into subintervals $K_1$, $L_1$, $K_2$, $L_2$, $K_3$, $O_1$, $O_2$, \ldots, $O_{2k - 1}$, $M_1$, $N_1$, $M_2$, $N_2$, $M_3$ of sizes $r$, $d - r$, $r$, $d - r$, $r$, $d$, $d$, \ldots, $d$, $r$, $d - r$, $r$, $d - r$, $r$, respectively.

Construct the pencils $\First(d - r, D'_1) \circeq L_1 \circeq N_1 \circeq L_2 \circeq N_2 \circeq \Last(d - r, D''_{2k})$ and $\Last(r, D'_1) \circeq K_2 \circeq M_2 \circeq \First(r, D''_{2k})$. This connects each vertex at position $i$ in $D'_1$ to the vertex at position $(i + r) \bmod d$ in $D''_{2k}$.

Construct also the pencils $\First(d - r, D'_2) \circeq \First(d - r, O_1)$, $\Last(r, D'_2) \circeq K_3 \circeq M_3 \circeq \Last(r, O_1)$, $D'_j \circeq O_{j - 1}$ for all $j$ with $3 \le j \le 2k$, $O_j \circeq D''_j$ for all $j$ with $1 \le j \le 2k - 2$, $\First(r, O_{2k - 1}) \circeq K_1 \circeq M_1 \circeq \First(r, D''_{2k - 1})$, and $\Last(d - r, O_{2k - 1}) \circeq \Last(d - r, D''_{2k - 1})$. This connects each vertex of $D'_j$ to the corresponding vertex of $D''_{j - 1}$ for all $j$ with $2 \le j \le 2k$.

The description of our tour is complete. Observe that by now all vertices of $I$ have degree two and that the connections between subintervals $D'_1$, $D'_2$, \ldots, $D'_{2k}$ and $D''_1$, $D''_2$, \ldots, $D''_{2k}$ are the same as the connections between subintervals $D_1$, $D_2$, \ldots, $D_{2k}$ and $D_2$, $D_3$, \ldots, $D_{2k + 1}$ in the proof of Lemma \ref{sde}.

(Or, in other words, the vertex at position $i'$ in $D'_{j'}$ and the vertex at position $i''$ in $D''_{j''}$ are connected by a path all of whose interior vertices are outside of $D'_1 \cup D'_2 \cup \cdots \cup D'_{2k} \cup D''_1 \cup D''_2 \cup \cdots \cup D''_{2k}$ if and only if an analogous claim holds for the vertex at position $i'$ in $D_{j'}$ and the vertex at position $i''$ in $D_{j'' + 1}$ in the proof of Lemma \ref{sde}. Though note also that in the latter setting the path will often be of zero length.)

From this point on, we establish that our construction does indeed yield a tour exactly as we did there.~\end{proof}

\medskip

\begin{figure}[ht] \centering \includegraphics[scale=0.45]{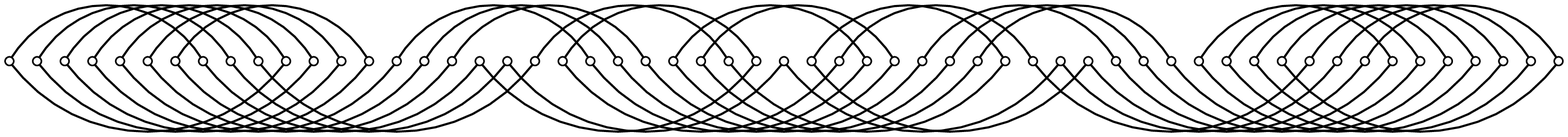} \caption{} \label{f:sdo} \end{figure}

Figure \ref{f:sdo} shows an example with $p = 7$ and $q = 10$, when $d = 3$, $r = 1$, and $\alpha = 3$.

The rest is not too difficult.

\medskip

\begin{lemma} \label{sd} The tours of projection graphs we constructed in Lemmas \ref{sde} and \ref{sdo} admit splits with widths as in Lemma \ref{ed}. \end{lemma} 

\medskip

\begin{proof} Breaking off either subgraph $\mathcal{B}'$ or $\mathcal{B}''$ as defined in the proofs of Lemmas \ref{sde} and \ref{sdo} yields a split with width $\alpha d$ when $\alpha$ is even and $(\alpha - 1)d$ when $\alpha$ is odd. Either way, we obtain a split whose width equals the greatest multiple of $2d$ which does not exceed $q$. Thus it satisfies the conditions of Lemma \ref{ed} in all cases except when $d = 1$ and $q$ is even. Then we break off a copy of $\mathcal{B}_{\alpha/2 - 1}$ instead, either on the left or on the right, for a split with width $q - 2d$. \end{proof}

\medskip

\begin{lemma} \label{bd} Suppose that $2p > q$. Suppose, furthermore, that both of $\Pi(p, q, n')$ and $\Pi(p, q, n'')$ are Hamiltonian. Then so is $\Pi(p, q, n' + n'')$. \end{lemma} 

\medskip

\begin{proof} Partition an integer interval $I$ of size $n' + n''$ into two subintervals $I'$ and $I''$ of sizes $n'$ and $n''$, respectively. Let $T'$ be a tour of $\Pi(p, q, I')$ and let $T''$ be a tour of $\Pi(p, q, I'')$.

Choose any vertex $u$ in $\Last(q - p, I')$. Then the pair of short edges $(u - p)$---$u$ and $(u - p + q)$---$(u + q)$ are forced in $T'$ and $T''$, respectively, as in the proof of Lemma \ref{se}. We delete them and we replace them with the pair of long edges $(u - p)$---$(u - p + q)$ and $u$---$(u + q)$ so that $T'$ and $T''$ merge together into a single tour of $\Pi(p, q, I)$, as in the proof of Lemma \ref{pmt}. \end{proof}

\medskip

\begin{lemma} \label{ptd} Suppose that $2p > q$. Then $\Pi(p, q, n)$ is Hamiltonian for all sufficiently large $n$. \end{lemma} 

\medskip

\begin{proof} Let $n^\star = \alpha(p + q) + d$ when $\alpha$ is even and $n^\star = \alpha(p + q) + 2d$ when $\alpha$ is odd. By Lemmas \ref{pmt}, \ref{es}, \ref{ed}, and \ref{sde}--\ref{bd} we obtain that $\Pi(p, q, n)$ is Hamiltonian for every $n$ of the form $n = n^\star + (p + q)i + 2qj$ where $i$ and $j$ are nonnegative integers. \end{proof}

\medskip

Combining Lemmas \ref{pkt}, \ref{pte}, and \ref{ptd}, we arrive at the following theorem.

\medskip

\begin{theorem} \label{pt} The projection graph $\Pi(p, q, n)$ is Hamiltonian for all sufficiently large positive integers~$n$. \end{theorem} 

\medskip

Let $\mu_\Pi$ be the smallest positive integer such that $\Pi(p, q, n)$ is Hamiltonian for all positive integers $n$ with $n \ge \mu_\Pi$. Conversely, the greatest positive integer $n$ such that $\Pi(p, q, n)$ is not Hamiltonian will then be $n = \mu_\Pi - 1$.

By the proof of Lemma \ref{pte}, when $2p < q$ we obtain that $\mu_\Pi \le (p + q) \cdot 2q < 3q^2$. Similarly, by the proof of Lemma \ref{ptd}, when $2p > q$ we obtain that $\mu_\Pi \le n^\star + (p + q) \cdot 2q < 6q^2$. Therefore, in all cases $\mu_\Pi \in \mathcal{O}(q^2)$.

We have learned enough about the tours of projection graphs for our immediate purposes. We will return to this topic once again in Section \ref{further}.

\section{Looms} \label{loom}

We define a \emph{loom} to be a pseudotour in a projection graph all of whose cycles alternate between short and long edges.

\begin{figure}[ht] \centering \includegraphics{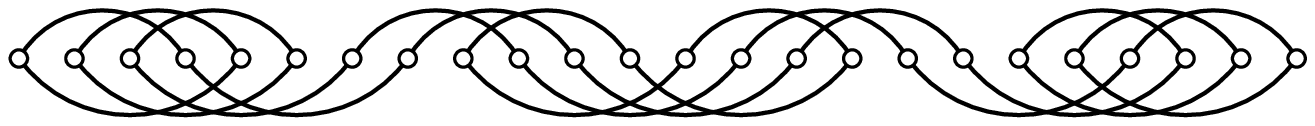} \caption{} \label{f:loom} \end{figure}

For example, Figure \ref{f:loom} shows a loom in $\Pi(3, 4, 24)$ with two cycles.

We can describe all looms explicitly.

\medskip

\begin{lemma} \label{2pq} The projection graph $\Pi(p, q, n)$ admits a loom if and only if $2pq$ divides $n$. Furthermore, for each such $n$ that loom is unique. \end{lemma}

\medskip

\begin{proof} The short edges of a loom form a perfect matching of $\Pi(p, q, n)$. So do the long ones. Conversely, the union of a perfect matching using only short edges and a perfect matching using only long edges is always a loom.

On the other hand, $\Pi(p, q, n)$ admits a perfect matching using only short edges if and only if $2p$ divides $n$. Furthermore, for each such $n$ that perfect matching is unique. The perfect matchings consisting entirely of long edges behave similarly. \end{proof}

\medskip

We denote the unique loom of $\Pi(p, q, 2kpq)$ by $\mathcal{L}_k$.

Partition the vertex set of $\Pi(p, q, 2kpq)$ into $2kq$ subintervals $P_1$, $P_2$, \ldots, $P_{2kq}$ of size $p$ each as well as into $2kp$ subintervals $Q_1$, $Q_2$, \ldots, $Q_{2kp}$ of size $q$ each. Then $\mathcal{L}_k$ is the union of all pencils of the form $P_{2i - 1} \circeq P_{2i}$ and all pencils of the form $Q_{2i - 1} \circeq Q_{2i}$.

In particular, $\mathcal{L}_k$ is the disjoint union of $k$ translation copies of $\mathcal{L}_1$ for all $k$.

One more explicit description of the edges of a loom will be useful to us in Section \ref{switch-ii}. Consider a loom $\mathcal{L}_k$ on vertices $[0; 2kpq - 1]$, and let $u$ be one of them. Then the short edge of $u$ in $\mathcal{L}_k$ joins it to $u + p$ when $\lfloor u/p \rfloor$ is even and to $u - p$ when $\lfloor u/p \rfloor$ is odd. Similarly, the long edge of $u$ in $\mathcal{L}_k$ joins it to $u + q$ when $\lfloor u/q \rfloor$ is even and to $u - q$ when $\lfloor u/q \rfloor$ is odd.

Since each cycle of a loom alternates between short and long edges, all such cycles are of even length. Therefore, a loom is always a bipartite graph.

The length of each cycle of a loom must also be a multiple of four because the alternate vertices along the cycle alternate between parities.

Let $\eta$ be the number of cycles in $\mathcal{L}_1$ and let $4\xi$ be the greatest length of a cycle in $\mathcal{L}_1$. Then the number of cycles in $\mathcal{L}_k$ is $k\eta$ and the greatest length of a cycle in $\mathcal{L}_k$ is $4\xi$ again. Furthermore, $\eta \le pq/2$ and $\xi \le pq/2$.

The point of looms will become clear in Section \ref{scarf}. Until then, we proceed to study some of the properties of cycles in projection graphs which alternate between short and long edges.

We call an oriented edge in a projection graph a \emph{move}. We say that move $u \to v$ points \emph{right} when $u < v$ and \emph{left} otherwise, when $u > v$.

We define a \emph{short bracket} to be a path in $\Pi(p, q, \mathbb{Z})$ of the form $u_1u_2u_3u_4$ such that both edges $u_1u_2$ and $u_3u_4$ are short, edge $u_2u_3$ is long, and moves $u_1 \to u_2$ and $u_3 \to u_4$ point in opposite directions. We define a \emph{long bracket} similarly, with ``short'' and ``long'' swapped.

\medskip

\begin{lemma} \label{lb} Every cycle in $\Pi(p, q, \mathbb{Z})$ which alternates between short and long edges contains a long bracket. \end{lemma} 

\medskip

\begin{proof} Let $C$ be our cycle and let $u_2$ be its leftmost vertex. Then the two neighbours of $u_2$ along $C$ will be $u_1 = u_2 + q$ and $u_3 = u_2 + p$. Furthermore, since $C$ alternates between short and long edges, the other neighbour of $u_3$ along $C$ will be $u_4 = u_3 + q$. \end{proof}

\medskip

Our proof of the analogous claim for short brackets is not quite so straightforward.

\medskip

\begin{lemma} \label{frog} Let $m$ and $n$ with $m < n$ be relatively prime positive integers and let $W$ be an oriented closed walk in $\Pi(m, n, \mathbb{Z})$ that contains only short moves to the right and long moves to the left. Furthermore, let $d$ be a positive integer with $d < m + n$. Then $W$ visits two vertices which differ by $d$. \end{lemma}

\medskip

The upper bound cannot be improved, as the appropriate orientation of the unique tour of $\Pi(m, n, m + n)$ serves as a counterexample for all $d$ with $d \ge m + n$.

(While working on this paper, the author submitted Lemma \ref{frog} as a problem proposal for the spring round of the 2020--2021 Tournament of Towns. It appeared there, as well as on the simultaneously held Moscow Mathematical Olympiad, before the paper was completed.)

\medskip

\begin{proof} Let $a$ be the smallest positive integer with $am \equiv d \pmod n$ and let $b = (am - d)/n$. Then $b$ is a nonnegative integer, $a \le n$, $b < m$, and $am - bn = d$.

Since $W$ is closed, it contains $kn$ short moves and $km$ long moves for some positive integer $k$. Consider all $k(m + n)$ subwalks of $W$ of length $a + b$, where we allow a subwalk to wrap around the coinciding beginning and ending of $W$.

By averaging over all such subwalks, we see that at least one of them contains at most $\lfloor \frac{m}{m + n} \cdot (a + b) \rfloor$ long moves and at least one of them contains at least $\lceil \frac{m}{m + n} \cdot (a + b) \rceil$ long moves.

On the other hand, as we pass from one subwalk to the next, deleting one move at the beginning and appending one move at the end, the number of long moves in our current subwalk changes by at most one. It follows that some subwalk of $W$ of length $a + b$ contains exactly $\lfloor \frac{m}{m + n} \cdot (a + b) \rfloor$ long moves.

However, $\lfloor \frac{m}{m + n} \cdot (a + b) \rfloor = b$ since the chain of inequalities $b < \frac{m}{m + n} \cdot (a + b) < b + 1$ simplifies to $0 < d < m + n$. Thus some subwalk of $W$ contains exactly $a$ short moves and $b$ long moves. Its endpoints will differ by $am - bn = d$. \end{proof}

\medskip

\begin{lemma} \label{sb} Every cycle in $\Pi(p, q, \mathbb{Z})$ which alternates between short and long edges contains a short bracket. \end{lemma} 

\medskip

\begin{proof} Consider any oriented closed walk $W$ in $\Pi(p, q, \mathbb{Z})$ which alternates between short and long moves and does not contain a short bracket. Then all short moves in $W$ must point in the same direction, without loss of generality to the right.

Let $u_1$, $u_2$, \ldots, $u_k$ be the midpoints of the short moves of $W$, in the order in which they occur along $W$. Then $u_{i + 1} - u_i$ equals $p + q$ when the long move between $u_i$ and $u_{i + 1}$ points to the right, and $p - q$ otherwise, when that long move points to the left. (With $i = k$ we allow the indices to wrap around cyclically, so that $u_{k + 1}$ is the same midpoint as $u_1$.)

Consequently, $u_1 \to u_2 \to \cdots \to u_k \to u_1$ is an oriented closed walk on the number line where every step is either a move of length $q - p$ to the left or a move of length $p + q$ to the right. It visits only integers when $p$ is even, and only half-integers otherwise, when $p$ is odd.

By Lemma \ref{frog} with $m = q - p$ and $n = p + q$, there exist indices $i$ and $j$ such that $|u_i - u_j| = p$. But then $(u_i + u_j)/2$ is an integer incident with two distinct short moves of $W$, and so also an integer which $W$ visits twice. Therefore, $W$ cannot be a simple cycle. \end{proof}

\medskip

Of course, both Lemmas \ref{lb} and \ref{sb} apply to all cycles in looms.

\section{Scarves} \label{scarf}

Consider a board $A = I \times J$ of size $m \times n$. Let $\Pi_X = \Pi(p, q, I)$ and $\Pi_Y = \Pi(p, q, J)$.

Take any edge $e = (x', y')$---$(x'', y'')$ of $L$ on $A$. Then $x'x''$ is an edge of $\Pi_X$, and we call it the \emph{$x$-projection} of $e$. Similarly, $y'y''$ is an edge of $\Pi_Y$ which we call the \emph{$y$-projection} of $e$.

Let $G_X$ be any subgraph of $\Pi_X$ and let $G_Y$ be any subgraph of $\Pi_Y$. Then we define the \emph{product} $G_X \bigcirc G_Y$ of $G_X$ and $G_Y$ to be subgraph of the leaper graph of $L$ on $A$ formed by all edges $e$ of $L$ on $A$ such that the $x$-projection of $e$ belongs to $G_X$ and the $y$-projection of $e$ belongs to $G_Y$.

\medskip

\begin{lemma} \label{ls} Let $\mathcal{L}$ be a loom in $\Pi_X$ and let $P$ be a pseudotour of $\Pi_Y$. Then $G = \mathcal{L} \bigcirc P$ is a pseudotour of $L$ on $A$. \end{lemma} 

\medskip

\begin{proof} Observe that an edge of $\mathcal{L}$ and an edge of $P$ are the two projections of an edge of $L$ on $A$ if and only if they have opposite types, so that one is short and the other one is long.

Take any cell $a = (x, y)$ of $A$. Let $x'$ and $x''$ be the neighbours of $x$ in $\mathcal{L}$ and let $y'$ and $y''$ be the neighbours of $y$ in $P$.

Since $\mathcal{L}$ is a loom, its edges $xx'$ and $xx''$ have opposite types. Consequently, exactly one of them combines with edge $yy'$ of $P$ to produce an edge of $G$ incident with $a$. The same reasoning applies to edge $yy''$ of $P$ as well. Therefore, every cell of $A$ has degree two in $G$. \end{proof}

\medskip

We call a pseudotour of $L$ on $A$ of this form a \emph{scarf}. We also call $P$ the scarf's \emph{pattern}.

Thus a loom produces a scarf according to a pattern.

\begin{figure}[t!] \centering \includegraphics{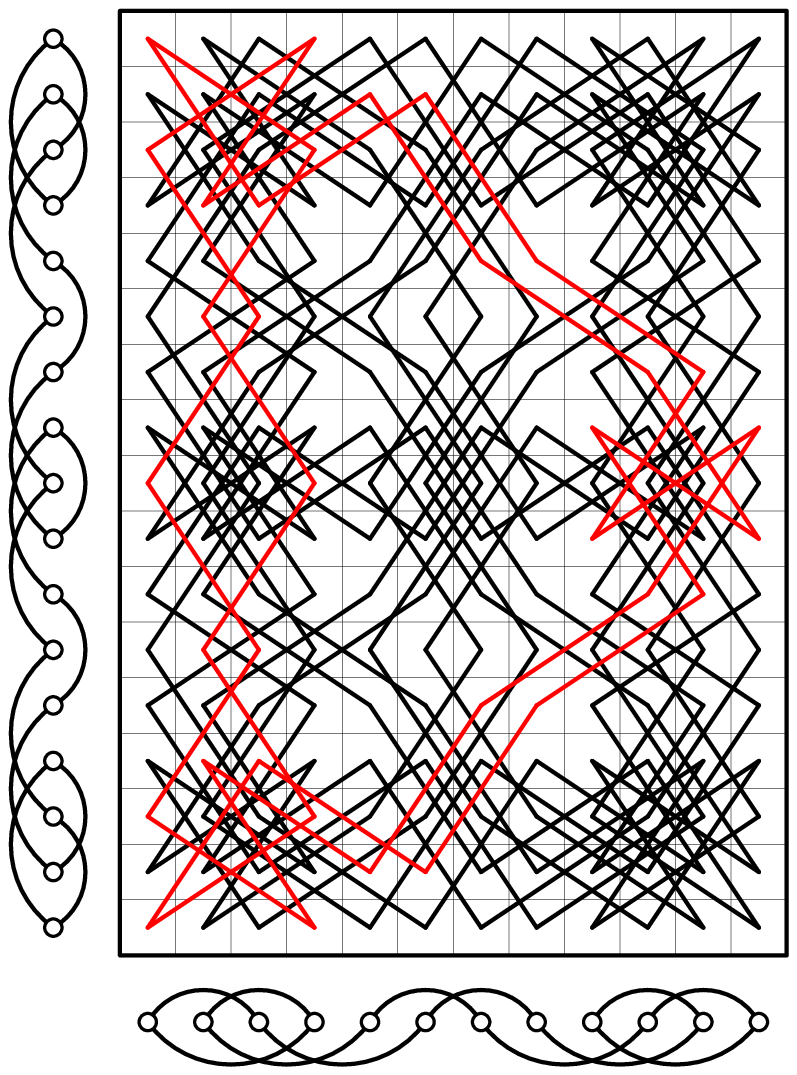} \caption{} \label{f:scarf} \end{figure}

For example, Figure \ref{f:scarf} shows a scarf of the zebra on the board of size $17 \times 12$, with the loom pictured below and the pattern on the left.

The idea of a scarf allows us to put together pseudotours of $L$ with great flexibility. In Sections \ref{switch-i}--\ref{leap}, we will learn how to modify some of these pseudotours so that they become tours.

Let $\mathcal{L}$, $P$, and $G = \mathcal{L} \bigcirc P$ be as in the statement of Lemma \ref{ls} and let $C = x_0x_1 \ldots x_{4\ell - 1}$ be a cycle in $\mathcal{L}$. (We already observed that the length of each cycle in a loom is a multiple of four in Section \ref{loom}.) Let $S'$ be the set of all cells of $A$ of the form $a = (x_i, y)$ such that $i$ and $a$ are of the same parity and let $S''$ be the set of all cells of $A$ of the form $a = (x_i, y)$ such that $i$ and $a$ are of opposite parities.

Then each cell of $S'$ is only joined in $G$ to other cells of $S'$. Consequently, $S'$ is the disjoint union of the vertex sets of several cycles of $G$. We say that these cycles form a \emph{braid}. Similarly, the cycles of $G$ making up $S''$ form a braid as well.

Thus we can group all cycles of $G$ into several braids so that exactly two braids are associated with each cycle of $\mathcal{L}$. Consequently, if $\mathcal{L} = \mathcal{L}_k$ then our scarf comprises $2k\eta$ braids altogether.

One corollary of this is that all visits to the same column by the same cycle of $G$ are on cells of the same parity. We proceed to examine visits to the same row as well.

\medskip

\begin{lemma} \label{odd} In the setting of Lemma \ref{ls}, suppose also that the height $m$ of $A$ is odd and the pattern $P$ is in fact a tour. Then each cycle of $G = \mathcal{L} \bigcirc P$ visits every row of $A$ exactly twice. Furthermore, the two visits are on cells of opposite parities. \end{lemma}

\medskip

Therefore, when $m$ is odd and $P$ is a tour each cycle of $G$ has length $2m$ and there are a total of $n/2$ such cycles. Furthermore, each braid associated with a cycle of $\mathcal{L}$ of length $4\ell$ contains exactly $\ell$ cycles. Thus, in particular, when additionally $p = 1$ each braid consists of a single cycle.

The scarf in Figure \ref{f:scarf} satisfies the conditions of Lemma \ref{odd}, and one of its cycles has been highlighted.

\medskip

\begin{proof} Let $(x_0, y_0)$ be any cell of $A$ and let $x_0x_1 \ldots x_{4\ell - 1}$ be the cycle of $\mathcal{L}$ which contains $x_0$. Let us trace out the cycle of $G$ which contains $a$ and let $(x_{i_0}, y_0) = (x_0, y_0)$, $(x_{i_1}, y_1)$, \ldots, $(x_{i_m}, y_m)$ be the cells of $A$ that we encounter over the course of the first $m$ steps of our journey.

By the proof of Lemma \ref{ls}, the two neighbours of $y_i$ in $P$ are $y_{i - 1}$ and $y_{i + 1}$ for all $i$ with $0 < i < m$. Therefore, every vertex of $P$ occurs exactly once among $y_0$, $y_1$, \ldots, $y_{m - 1}$ and, furthermore, $y_m = y_0$.

Since $4\ell$ is even and $m$ is odd, $i_m$ is odd as well. In particular, $i_m \neq i_0$ and so $(x_{i_m}, y_m) \neq (x_0, y_0)$.

For all $i$, define $\overline{i} = (i_m - i) \bmod 4\ell$. Observe that, since $i_m$ is odd, edges $x_ix_{i + 1}$ and $x_{\overline{i}}x_{\overline{i + 1}}$ of $\mathcal{L}$ are always of the same type, so that either both of them are short or both of them are long. (With $i = 4\ell - 1$ we allow the indices to wrap around cyclically, so that $x_{4\ell}$ is the same vertex of $\mathcal{L}$ as $x_0$.)

Therefore, over the course of our next $m$ steps along the cycle of $a$ in $G$ beyond cell $(x_{i_m}, y_m) = (x_{\overline{i_0}}, y_0)$ we will successively visit cells $(x_{\overline{i_1}}, y_1)$, $(x_{\overline{i_2}}, y_2)$, \ldots, $(x_{\overline{i_m}}, y_m) = (x_0, y_0)$. Or, in other words, on step $2m$ we will return back to cell $(x_0, y_0)$ and close the cycle.

Since $P$ is a tour of $\Pi_Y$, our route does visit every row of $A$ exactly twice. Furthermore, since the two visits are always separated by $m$ moves of $L$ and $m$ is odd, for each row they do indeed occur on cells of opposite parities. \end{proof}

\section{Switches I} \label{switch-i}

Let $G$ be a pseudotour of $L$ on $A$.

Consider a cycle $S$ of $L$ on $A$ such that $S$ alternates between edges in $G$ and edges outside of $G$. Suppose, furthermore, that all edges of $S$ in $G$ belong to pairwise distinct cycles of $G$. Then we call $S$ a \emph{switch}.

Let us modify $G$ as follows. First we delete all edges of $S$ in $G$. Then we replace them with the edges of $S$ outside of $G$. The net result is that all cycles of $G$ touched by $S$ become stitched together into one single longer cycle. This is what we call \emph{flipping} a switch.

We might hope to flip several switches one by one so that all cycles of $G$ merge together and we obtain a tour. This method of construction is well-known in the context of knight tours. For details see \cite{J2}, where the term \emph{linkage polygon} is used instead of ``switch''.

The smallest possible length of a switch is four, when it touches only two cycles of $G$. Then we call $S$ a \emph{rhombus switch} because the centers of its cells are the vertices of a rhombus.

(Geometrically, there are three distinct types of such rhombuses: (a) When two opposite vertices are in the same row and the other two are in the same column; (b) When two opposite vertices are in the same diagonal of slope $1$ and the other two are in the same diagonal of slope $-1$; and (c) When the four vertices form a slanted square. Of these, types (b) and (c) do occur in scarves whereas type (a) does not. However, type (a) can occur in other pseudotours.)

Let $\mathcal{S} = \{S_1$, $S_2$, \ldots, $S_k\}$ be a system of switches and let $\mathcal{C} = \{C_1$, $C_2$, \ldots, $C_\ell\}$ be the set of all cycles of $G$ that these switches touch. We say that $\mathcal{S}$ \emph{connects} $\mathcal{C}$ when the union of all switches of $\mathcal{S}$ and all cycles of $\mathcal{C}$ is a connected graph. In the special case when $\mathcal{C}$ consists of all cycles of $G$, we say simply that $\mathcal{S}$ connects $G$.

Clearly, if we can flip several switches of $\mathcal{S}$ in succession so that all cycles of $\mathcal{C}$ become stitched together then $\mathcal{S}$ must necessarily connect $\mathcal{C}$. The converse, however, is false in general because the first few flips might transform $G$ in such a way that some of the cycles of $\mathcal{S}$ no longer satisfy the definition of a switch in the transformed pseudotour.

But it does hold in some special cases. We find the next lemma in \cite{B}. (Though its wording and proof are different there.)

\medskip

\begin{lemma*}{\textbf{R}} (Rhombus switches.) Let $\mathcal{S}$ be a connecting system of pairwise edge-disjoint rhombus switches. Then we can flip several switches of $\mathcal{S}$ in succession so as to stitch all cycles touched by the switches of $\mathcal{S}$ together into one single longer cycle. \end{lemma*}

\medskip

\begin{proof} Let $F$ be the graph whose vertices are all cycles touched by the switches of $\mathcal{S}$ and where two such cycles are joined by an edge if and only if some switch of $\mathcal{S}$ touches both of them. Choose a spanning tree $T$ of $F$ and, for each edge of $T$, choose one switch of $\mathcal{S}$ which touches the two cycles joined by that edge. Then flip all of these switches one by one. \end{proof}

\medskip

We go on to construct connecting systems of rhombus switches for the braids of certain scarves.

When $p = 1$ and our pattern is a tour of odd length, every braid consists of a single cycle as we found in Section \ref{scarf}. Suppose, throughout the rest of this section, that $p \ge 2$.

\begin{figure}[ht] \centering \includegraphics{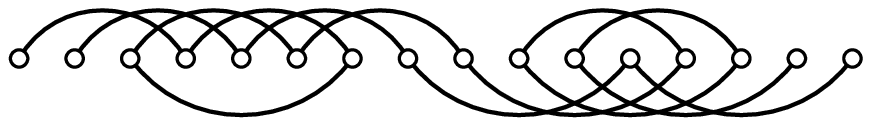} \caption{} \label{f:ea} \end{figure} 

We define the extension $\mathcal{E}$ as follows. When $L$ is distinct from the antelope, $\mathcal{E}$ is as constructed in Lemmas \ref{ee} and \ref{ed}. Otherwise, in the special case of the antelope, we use the extension of width two and length fourteen consisting of the two paths $0$---$3$---$6$---$2$---$5$---$8$---$12$---$9$---$13$---$10$---$14$ and $1$---$4$---$7$---$11$---$15$, as shown in Figure \ref{f:ea}.

\medskip

\begin{lemma} \label{swb} Suppose that $p \ge 2$. Let $\mathcal{L}$ be a loom and let $P$ be a tour in a projection graph of odd size which contains a chain $H$ of sufficiently many copies of $\mathcal{E}$. Consider the scarf $G = \mathcal{L} \bigcirc P$. Then for every braid $\mathcal{C}$ of $G$ there exists a system $\mathcal{S}$ of pairwise edge-disjoint rhombus switches for $G$ such that $\mathcal{S}$ connects $\mathcal{C}$ and the $y$-projections of the switches of $\mathcal{S}$ are all contained within $H$. \end{lemma} 

\medskip

\begin{proof} We designate two long edges $e'$ and $e''$ of $\mathcal{E}$ as follows.

\smallskip

\emph{Case 1}. $2p < q$. In the notation of Lemma \ref{ee} and its proof, find a vertex $u$ in $\First(p, J')$ for which the smallest positive integer $i$ with $u + pi \in I_\text{Mid}$ is odd. The existence of such a vertex is guaranteed by the fact that $p$ does not divide $q$. Then denote $v = u + pi$ and set $e' = (v - q)$---$v$ and $e'' = (v + p - q)$---$(v + p)$.

Observe that the path of $e'$ in $\mathcal{E}$ consists of two blocks of an odd number of short edges each together with two long edges, whereas the path of $e''$ in $\mathcal{E}$ consists simply of two long edges.

\smallskip

\emph{Case 2}. $2p > q$ and $L$ is distinct from the antelope. In the notation of Lemma \ref{ed} and its proof, set $w = d$ when $s = q - 2d$ and $w = \max\{0, s - p\}$ otherwise, when $s > q - 2d$. Then also set $e' = (w + p)$---$(w + p + q)$ and $e'' = w$---$(w + q)$.

Observe that our choice of $w$ ensures that the path of $e'$ in $\mathcal{E}$ is of length four and alternates between short and long edges, whereas the path of $e''$ in $\mathcal{E}$ consists simply of two long edges.

\smallskip

\emph{Case 3}. $L$ is the antelope. Then set $e' = 10$---$14$ and $e'' = 7$---$11$.

\smallskip

Let $\mathcal{E}_0$, $\mathcal{E}_1$, \ldots, $\mathcal{E}_{k - 1}$ be the copies of $\mathcal{E}$ forming $H$, in order from left to right, and let $e'_i = y'_iz'_i$ with $y'_i < z'_i$ and $e''_i = y''_iz''_i$ with $y''_i < z''_i$ be the copies of $e'$ and $e''$, respectively, in $\mathcal{E}_i$.

Let $C = x_0x_1 \ldots x_{4\ell - 1}$ be a cycle of $\mathcal{L}$. By Lemma \ref{sb}, $C$ contains a short bracket. Suppose, without loss of generality, that $x_0x_1x_2x_3$ is one such. Consequently, in particular all edges of $C$ of the form $x_{2i}x_{2i + 1}$ are short and all of its other edges are long.

Let $\mathcal{C}$ be one of the two braids of $C$. Then $C$ contains either cell $(x_0, y'_0)$ or cell $(x_2, y'_0)$. We consider the former case in detail, and the latter case is similar.

By Lemma \ref{odd}, the cells $(x_{4i}, y'_0)$ belong one each to the cycles of $\mathcal{C}$. Let $\mathcal{C} = \{C_0, C_1, \ldots, C_{\ell - 1}\}$ so that cell $(x_{4i}, y'_0)$ belongs to cycle $C_i$. We allow the indices of both the vertices of $C$ and the cycles of $\mathcal{C}$ to wrap around cyclically, so that for example $x_{-1}$ is the same vertex as $x_{4\ell - 1}$ and $C_\ell$ is the same cycle as $C_0$.

We proceed to show that cell $(x_{4i}, y'_j)$ belongs to cycle $C_{i - j}$ for all $i$ and $j$. It suffices to verify that cells $(x_{4i}, y'_j)$ and $(x_{4(i + 1)}, y'_{j + 1})$ belong to the same cycle of $G$ for all $i$ and all $j$ with $0 \le j \le k - 2$. We carry out the verification in detail for Case 1 of the definition of $e'$ and $e''$, when $2p < q$, and the other two cases are fully analogous.

Consider the path of $H$ which contains $e'_j$ and $e'_{j + 1}$. Let the portion of it which extends from $y'_j$ to $y'_{j + 1}$ be $V = v_0v_1 \ldots v_{4s}$ so that $v_0 = y'_j$, $v_1 = z'_j$, $v_1v_2 \ldots v_{2s}$ is a block of an odd number of short edges, $v_{2s}v_{2s + 1}$ is a long edge, $v_{2s + 1}v_{2s + 2} \ldots v_{4s}$ is another block of an odd number of short edges, and, finally, $v_{4s} = y'_{j + 1}$.

Consider also the walk $U = u_0u_1 \ldots u_{4s}$ in $C$ which starts from $x_{4i}$, goes to $x_{4i + 1}$, oscillates between $x_{4i + 1}$ and $x_{4i + 2}$ for $2s - 1$ steps, continues on to $x_{4i + 3}$, oscillates once again between $x_{4i + 3}$ and $x_{4(i + 1)}$ for $2s - 1$ more steps, and eventually settles at $x_{4(i + 1)}$. Explicitly, $u_0 = x_{4i}$, $u_1 = u_3 = \cdots = u_{2s - 1} = x_{4i + 1}$, $u_2 = u_4 = \cdots = u_{2s} = x_{4i + 2}$, $u_{2s + 1} = u_{2s + 3} = \cdots = u_{4s - 1} = x_{4i + 3}$, and $u_{2s + 2} = u_{2s + 4} = \cdots = u_{4s} = x_{4(i + 1)}$.

Then the corresponding edges of $U$ and $V$ are always of opposite types, so that one is short and the other one is long. Consequently, $(u_0, v_0)$---$(u_1, v_1)$---$\cdots$---$(u_{4s}, v_{4s})$ is a path in $G$ which connects cells $(u_0, v_0) = (x_{4i}, y'_j)$ and $(u_{4s}, v_{4s}) = (x_{4(i + 1)}, y'_{j + 1})$. Therefore, these cells do indeed belong to the same cycle of $G$.

By Lemma \ref{odd}, the cells $(x_{4i + 2}, z''_0)$ also belong one each to the cycles of $\mathcal{C}$. Let $\mathcal{C} = \{D_0, D_1, \ldots, D_{k - 1}\}$ so that cell $(x_{4i + 2}, z''_0)$ belongs to cycle $D_i$.

By the same method of analysis as above, we obtain that cell $(x_{4i + 2}, z''_j)$ belongs to cycle $D_i$ of $\mathcal{C}$ for all $i$ and $j$.

Consider the cycles $S_j = (x_0, y'_j)$---$(x_1, z'_j)$---$(x_2, z''_j)$---$(x_3, y''_j)$---$(x_0, y'_j)$ of $L$ on $A$. We just established that two opposite edges of $S_j$ belong to cycles $C_{-j}$ and $D_0$ of $\mathcal{C}$. Therefore, $S_j$ is a rhombus switch for all $j$ such that $C_{-j}$ and $D_0$ are distinct. Furthermore, provided that $k \ge \xi$, where $\xi$ is as in Section \ref{loom}, among $S_0$, $S_1$, \ldots, $S_{k - 1}$ there will be rhombus switches connecting $D_0$ to all other cycles of $\mathcal{C}$. \end{proof}

\medskip

By now we know how to construct connecting system of switches within braids. We devote the next section to the construction of connecting systems of switches which bring together different braids.

\section{Switches II} \label{switch-ii}

Let $u$ be any vertex of a pseudotour in a projection graph. We say that $u$ is a \emph{straight} when its two edges in the pseudotour lie on opposite sides of it, and we say that $u$ is a \emph{turn} otherwise, when they lie on the same side.

Consider the loom $\mathcal{L}_2$ on vertices $I = [0; 4pq - 1]$. By our explicit description of the edges of a loom from Section \ref{loom}, $u$ is a straight in $\mathcal{L}_2$ when $\tau(u) = \lfloor u/p \rfloor + \lfloor u/q \rfloor$ is odd and a turn otherwise, when $\tau(u)$ is even.

Partition $I$ into two halves $I'$ and $I''$, each of size $2pq$. Then each one of $I'$ and $I''$ contains equally many straights and turns since $\tau(u + pq) = \tau(u) + p + q$ and so each pair of vertices that differ by $pq$ consists of one straight and one turn. Furthermore, two corresponding vertices of $I'$ and $I''$ will always be of the same type because $\tau(u + 2pq) = \tau(u) + 2(p + q)$.

Let $u$ and $u + 2pq$ be two corresponding straights in $I'$ and $I''$, respectively. Consider the cycle $u$---$(u + 2p)$---$(u + 4p)$---$\cdots$---$(u + 2pq)$---$(u + 2(p - 1)q)$---$(u + 2(p - 2)q)$---$\cdots$---$u$ in $\Pi(2p, 2q, I)$. We call a cycle of this form a \emph{shaft}. (That a shaft is indeed a simple cycle, rather than a closed walk which revisits vertices, follows because $p$ and $q$ are relatively prime.) We define the parity of a shaft to be the parity of all of its vertices.

Let $P$ be a pseudotour in the projection graph $\Pi(p, q, J)$ and let $G = \mathcal{L}_2 \bigcirc P$ be the scarf of loom $\mathcal{L}_2$ and pattern $P$ on the board $A = I \times J$.

We say that a vertex of $P$ is \emph{pure} when it is incident with two edges of $P$ of the same type, and we say that it is \emph{mixed} otherwise, when it is incident with one short and one long edge of $P$. (Thus in particular a vertex of $P$ cannot be simultaneously pure and a turn.)

Let $v$ be a mixed vertex of $P$ with short-edge neighbour $v'$ and long-edge neighbour $v''$ and let $S$ be a shaft. To each edge $u'u''$ of $S$, we assign a path of $L$ on $A$ as follows. Let $u_\text{Mid} = (u' + u'')/2$. When $u'u''$ is short, so that $|u' - u''| = 2p$, we assign to it the path $(u', v)$---$(u_\text{Mid}, v'')$---$(u'', v)$. Otherwise, when $u'u''$ is long, so that $|u' - u''| = 2q$, we assign to it the path $(u', v)$---$(u_\text{Mid}, v')$---$(u'', v)$.

The concatenation of all of these paths will then be a cycle of $L$ on $A$ which we call the \emph{comb} of index $v$ and shaft $S$.

(We think of each path assigned to an edge of $S$ as a ``tooth'' in our comb.)

Let $u$ and $u + 2pq$ be the two outermost vertices of $S$. Then the comb with index $v$ and shaft $S$ is explicitly $(u, v)$---$(u + q, v')$---$(u + 2q, v)$---$(u + 3q, v')$---$(u + 4q, v)$---$\cdots$---$(u + 2pq, v)$---$(u + (2q - 1)p, v'')$---$(u + 2(q - 1)p, v)$---$(u + (2q - 3)p, v'')$---$(u + 2(q - 2)p, v)$---$\cdots$---$(u, v)$.

Let $m$ be the size of $J$.

\medskip

\begin{lemma} \label{comb} Suppose that $m$ is odd and the pattern $P$ is in fact a tour. Then every comb is also a switch. \end{lemma}

\medskip

\begin{proof} That a comb alternates between edges in $G$ and edges outside of $G$ is clear. That the edges of a comb in $G$ belong to pairwise distinct cycles of $G$ then follows by Lemma \ref{odd}. \end{proof}

\medskip

\begin{figure}[t!] \centering \includegraphics{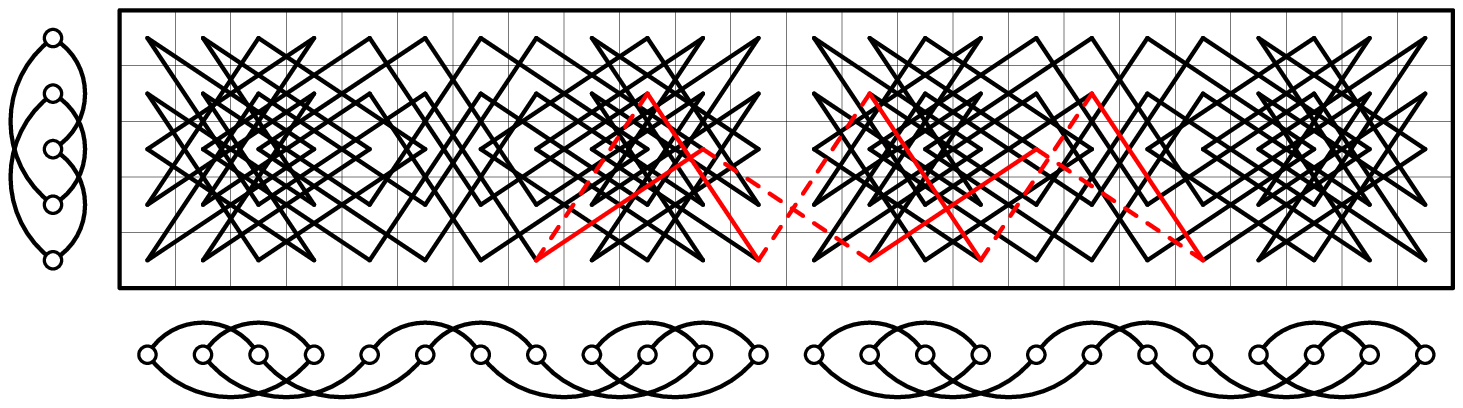} \caption{} \label{f:comb} \end{figure}

For example, Figure \ref{f:comb} shows the comb with index $0$ and shaft $7$---$11$---$15$---$19$---$13$---$7$ in the scarf of the zebra that we obtain when $P$ is the entire projection graph $\Pi(2, 3, 5)$.

We proceed to construct a connecting system of comb switches which touches at least one cycle out of every braid.

\medskip

\begin{lemma} \label{us} The union of all even shafts is a connected subgraph of $\Pi(2p, 2q, I)$. \end{lemma} 

\medskip

\begin{proof} Let $u_1 < u_2 < \cdots < u_{pq/2}$ be all even straights in $I'$ and let $S_i$ be the shaft with outermost vertices $u_i$ and $u_i + 2pq$. It suffices to prove that $S_i$ and $S_{i + 1}$ have a common vertex for all $i$ with $1 \le i < pq/2$.

Fix $i$ and let $u$ be the smallest positive integer such that $u > u_i$, $u \equiv u_i \pmod{2p}$, and $u \equiv u_{i + 1} \pmod{2q}$. Then $u$ is a common vertex of $S_i$ and $S_{i + 1}$ unless $u < u_{i + 1}$, in which case $u_{i + 1} - u_i = (u - u_i) + (u_{i + 1} - u) \ge 2p + 2q$. We proceed to show that in fact always $u_{i + 1} - u_i \le 2p + 2$, and so the latter inequality cannot occur.

When $p = 1$, explicitly $u_{i + 1} - u_i = 2$ for all $i$ with $1 \le i < pq/2$. Thus suppose, from this point on throughout the rest of the proof, that $p \ge 2$. Suppose also, for the sake of contradiction, that $u_i + 2j$ is a turn in $I'$ for all $j$ with $1 \le j \le p + 1$.

Since the integer interval $[u_i + 3; u_i + 2p + 2]$ is of size $2p$, it contains exactly two multiples of $p$, say $u'$ and $u''$. Let $v'$ be whichever one of $u'$ and $u' + 1$ is even, and define $v''$ similarly. Then both of $v'$ and $v''$ are among $u_i + 4$, $u_i + 6$, \ldots, $u_i + 2p + 2$, and they satisfy $v' \bmod p \in \{0, 1\}$ and $v'' \bmod p \in \{0, 1\}$.

Since $p \ge 2$, consequently $\lfloor (v' - 2)/p \rfloor = \lfloor v'/p \rfloor - 1$. On the other hand, since both of $v' - 2$ and $v'$ are turns by assumption, $\tau(v' - 2)$ and $\tau(v')$ must have the same parity. Therefore, $\lfloor (v' - 2)/q \rfloor = \lfloor v'/q \rfloor - 1$ and so $v' \bmod q \in \{0, 1\}$ as well.

Observe that $v' \neq 0$ and also $v' \neq pq$ because $\tau(pq) = p + q$ is odd. Thus we can rule out the case of $v' \bmod p = v' \bmod q = 0$. On the other hand, since $v'$ is even, it must necessarily be congruent to zero modulo whichever one of $p$ and $q$ is even. Consequently, it must also be congruent to one modulo whichever one of $p$ and $q$ is odd.

The same reasoning applies to $v''$ as well. Therefore, $v' \equiv v'' \pmod p$ and $v' \equiv v'' \pmod q$. However, $|v' - v''| \le (u_i + 2p + 2) - (u_i + 4) = 2p - 2 < pq$. We have arrived at a contradiction. \end{proof}

\medskip

\begin{lemma} \label{ms} Strictly more than half of all even vertices of $\mathcal{L}_2$ belong to even shafts. \end{lemma} 

\medskip

\begin{proof} We already know that exactly half of all even vertices of $\mathcal{L}_2$ are straights, and consequently belong to even shafts. To establish a strict inequality, it suffices to exhibit an even turn in a shaft. When $p$ is even, we consider vertex $2pq - p$ which belongs to the shaft with outermost vertices $p$ and $p + 2pq$. Otherwise, when $p$ is odd, we consider vertex $2pq + p - 1$ which belongs to the shaft with outermost vertices $2pq - p - 1$ and $4pq - p - 1$. \end{proof}

\medskip

Of course, both Lemmas \ref{us} and \ref{ms} apply to odd shafts as well.

\medskip

\begin{lemma} \label{sbb} Suppose that $m$ is odd and the pattern $P$ is in fact a tour. Let $v_0$ be an even mixed vertex of $P$ and let $v_1$ be an odd mixed vertex of $P$. Consider the system $\mathcal{S}$ of all comb switches whose index is one of $v_0$ and $v_1$. Then $\mathcal{S}$ is a connecting system whose switches touch at least one cycle out of every braid. \end{lemma} 

\medskip

\begin{proof} For all parities $i$ and $j$, let $\mathcal{S}_{ij}$ be the subsystem of $\mathcal{S}$ consisting of all comb switches in $\mathcal{S}$ with index of parity $i$ and shaft of parity $j$.

By Lemma \ref{us}, the union of all switches of $\mathcal{S}_{ij}$ is a connected subgraph of the leaper graph of $L$ on $A$. Therefore, $\mathcal{S}_{ij}$ is a connecting system of switches for all $i$ and $j$.

By Lemma \ref{ms}, strictly more than half of all cells in row $v_i$ with parity $i + j$ belong to the switches of $\mathcal{S}_{ij}$. By Lemma \ref{odd}, it follows that these switches touch strictly more than half of all cycles of $G$. Consequently, for each two out of $\mathcal{S}_{00}$, $\mathcal{S}_{01}$, $\mathcal{S}_{10}$, and $\mathcal{S}_{11}$, there exists a cycle of $G$ that the switches of both of them touch. Therefore, $\mathcal{S}$ is a connecting system of switches as well.

Finally, let $C$ be any cycle of $\mathcal{L}_2$. By either Lemma \ref{lb} or \ref{sb}, there exists a straight $u$ in $C$. Then cells $(u, v_0)$ and $(u, v_1)$ belong one each to the two braids of $C$ and the switches of $\mathcal{S}$ touch the cycles of $G$ which contain them. \end{proof}

\medskip

For our purposes in the next section, it will be more convenient to abstract away some of the content of Lemma \ref{sbb} so that we can view things in terms of braids rather than in terms of individual cycles.

We call each subset of $I$ formed by the alternate vertices of some cycle of $\mathcal{L}_2$ a \emph{signature}. To each braid $\mathcal{C}$ of $G$, we assign the signature which consists of all elements $u$ of $I$ such that $\mathcal{C}$ contains the even cells of column $u$.

Let $C$ be a comb switch for $G$ with index $v$ of parity $d$ and shaft $S$. Consider the set $\psi(C)$ of the signatures of all braids $\mathcal{C}$ such that $C$ touches at least one cycle of $\mathcal{C}$. Observe that $\psi(C)$ depends only on $d$ and $S$. It does not depend on the specific value of $v$, or the pattern $P$ of our scarf, or the height $m$ of the board. Thus we can also write simply $\psi(d, S)$.

Let $\mathcal{H}$ be the hypergraph whose vertices are all signatures and whose hyperedges are all sets of the form $\psi(d, S)$ over both parities $d$ and all shafts $S$.

\medskip

\begin{lemma} \label{ch} The hypergraph $\mathcal{H}$ is connected. \end{lemma} 

\medskip

\begin{proof} This is a direct corollary of Lemma \ref{sbb}. \end{proof}

\medskip

One might wonder why we bother to introduce comb switches when rhombus switches are so much more convenient to work with. The reason is that for many leapers there does not exist a scarf connected by a system consisting entirely of rhombus switches.

We can verify this as follows. Let $\mathcal{L}^\star_k$ be the union of $\mathcal{L}_k$ with all cycles $D$ in $\Pi(p, q, 2kpq)$ such that $D$ is of length four, $D$ alternates between short and long edges, and two opposite edges of $D$ are in $\mathcal{L}_k$. Suppose that a system of rhombus switches connects some scarf on loom $\mathcal{L}_k$. Then the $x$-projections of the switches in that system would ensure that $\mathcal{L}^\star_k$ is connected. However, in fact $\mathcal{L}^\star_k$ is disconnected for all $k$ for all $95$ skew free leapers $L$ with $p \ge 3$ and $p + q \le 35$.

One might also wonder why we base our construction on the loom $\mathcal{L}_2$ rather than the smallest loom $\mathcal{L}_1$. This time around the answer is simpler: By the nonexistence results of \cite{Kn} and \cite{O} we cited in the introduction it follows that, when $L$ belongs to either one of the two exceptional families $p = 1$ and $q = p + 1$, it can never tour a board of width~$2pq$.

\section{Tours in Leaper Graphs} \label{leap}

We are ready to tackle Theorem \ref{4pq}.

The main difficulty is as follows. Over the course of Sections \ref{switch-i} and \ref{switch-ii}, we constructed two systems of switches whose union connects some scarves. However, the rhombus switches lemma does not apply because one system consists of comb switches, and so we cannot use that union directly in the way outlined at the beginning of Section \ref{switch-i}.

We get around this obstacle by constructing one big scarf out of many smaller subscarves and then flipping at most one comb switch within each one of them so as to avoid unwelcome interference. This allows the comb switches to do their part of the job without getting in each other's way. The rest we handle with the help of the rhombus switches lemma.

\medskip

\begin{proof*}{Proof of Theorem \ref{4pq}} Let $m$ be a sufficiently large positive integer. We will show that $L$ tours the board of size $m \times 4pq$.

Let $I = [0; 4pq - 1]$ and let $\mathcal{L}_2$ be the loom on vertex set $I$.

Since $m$ is sufficiently large, we can write $m = m_1 + m_2 + \cdots + m_k$ so that $k$ and all of the $m_i$ are sufficiently large as well and all of the $m_i$ are odd. Consider an integer interval $J$ of size $m$, and partition it into subintervals $J_1$, $J_2$, \ldots, $J_k$ of sizes $m_1$, $m_2$, \ldots, $m_k$, respectively.

Let $P_i$ be a tour of $\Pi(p, q, J_i)$ constructed as in Sections \ref{proj-i} and \ref{proj-ii}. Let also $G_i = \mathcal{L}_2 \bigcirc P_i$ be the scarf of loom $\mathcal{L}_2$ and pattern $P_i$ on the board $A_i = I \times J_i$.

Let $P$ be the disjoint union of the $P_i$ and let $G$ be the disjoint union of the $G_i$. Then $P$ is a pseudotour of $\Pi(p, q, J)$ and $G = \mathcal{L}_2 \bigcirc P$ is the scarf of loom $\mathcal{L}_2$ and pattern $P$ on the board $A = I \times J$ of size $m \times 4pq$.

The case of the knight is exceptional, and so we deal with it first. Let $J_i = [a_i; b_i]$. For each $i$, flip the rhombus switch with opposite vertices $(2, a_i)$ and $(5, a_i + 3)$ as well as the comb switch with index $a_i$ and shaft $1$---$3$---$5$---$1$. This transforms each subscarf $G_i$ into a tour of the subboard $A_i$. Then for each $i$ with $1 \le i < k$ also flip the rhombus switch with opposite vertices $(1, b_i)$ and $(2, a_{i + 1})$. This causes all such tours to merge together into a single longer tour.

Suppose, from this point on throughout the rest of the proof, that $L$ is distinct from the knight.

Let $\mathcal{E}$ be the extension of Lemma \ref{ee} when $p = 1$ and the extension of Lemma \ref{swb} when $p \ge 2$. We impose a series of additional constraints on the patterns $P_i$ for all $i$, as follows: (a) $P_i$ contains a chain $H_i$ of sufficiently many copies of $\mathcal{E}$; (b) $P_i$ contains one more instance $E_i$ of $\mathcal{E}$ outside of the chain $H_i$; and (c) When $2p < q$, furthermore $P_i$ contains one more instance $F_i$ of $\mathcal{E}$ on vertex set $\First(2(p + q), J_i)$ distinct from $E_i$ and outside of the chain $H_i$.

Clearly, we can indeed find a pattern $P_i$ satisfying all three constraints (a), (b), and (c) whenever $m_i$ is sufficiently large.

We proceed to construct three separate systems of switches for $G$.

For the first system, we begin with the case when $p \ge 2$. Let $\mathcal{H}$ be the hypergraph of Lemma \ref{ch} and let $e_1$, $e_2$, \ldots, $e_\ell$ be the hyperedges of a minimal connected spanning subgraph of $\mathcal{H}$. (Thus $\ell \le 4\eta - 1$, where $\eta$ is as in Section \ref{loom}.) Since $k$ is sufficiently large, we can safely assume that $k \ge \ell$.

Observe that in this case $\mathcal{E}$ contains two mixed vertices of opposite parities. (We defined pure and mixed vertices in Section \ref{switch-ii}.) Thus for each $i$ with $1 \le i \le \ell$ we can find a comb switch $C_i$ in the subscarf $G_i$ such that $\psi(C_i) = e_i$ and the $y$-projection of $C_i$ is contained entirely within $E_i$. Let $\mathcal{S}_\text{I}$ be the system of all such comb switches.

Otherwise, when $p = 1$, all mixed vertices of $\mathcal{E}$ are even and so we alter the construction accordingly. Choose one mixed vertex $v_i$ of $P_i$ both of whose edges are in $E_i$ for $i = 1$ and $i = 2$. Let $C_i$ be the comb switch of index $v_i$ and shaft with opposite vertices $1$ and $2q + 1$ when $i = 1$ as well as $2q - 2$ and $4q - 2$ when $i = 2$. Let also $\ell = 2$ and $\mathcal{S}_\text{I} = \{C_1, C_2\}$. Then the hyperedges $e_1$ and $e_2$ which correspond to $C_1$ and $C_2$ will form a minimal connected spanning subgraph of $\mathcal{H}$ regardless of the parities of $v_1$ and $v_2$.

For the second system, once again we begin with the case when $p \ge 2$. By Lemma \ref{swb}, for each $i$ and each braid of $G_i$ there exists a system of pairwise edge-disjoint rhombus switches which connects that braid and the $y$-projections of all of whose switches are contained entirely within $H_i$. Let $\mathcal{S}_\text{II}$ be the union of all such systems over all $i$ and all braids of $G_i$. Then $\mathcal{S}_\text{II}$ is also a system of pairwise edge-disjoint rhombus switches.

Otherwise, when $p = 1$, by our observations in Section \ref{scarf} each braid of each subscarf $G_i$ consists of a single cycle and so we let $\mathcal{S}_\text{II} = \varnothing$. Thus in this case the ``sufficiently many'' part of constraint (a) could safely mean ``zero''.

For the third system, let $w = q - 1$ when $3p < q$ and $w = \lfloor (p + q)/2 \rfloor$ when $3p > q$. For all $i$ with $1 \le i < k$, let $w_{1, i}$ be the vertex at position $w$ in $J_{i + 1}$ and define also $w_{2, i} = w_{1, i} - p$, $w_{3, i} = w_{1, i} - p - q$, and $w_{4, i} = w_{1, i} - q$. Then $w_{1, i}w_{2, i}$ is a short edge of $P_{i + 1}$ and $w_{3, i}w_{4, i}$ is a short edge of $P_i$. This follows by constraint (c) when $3p < q$ and because the two edges are forced in the two patterns otherwise, when $3p > q$.

(The reason we impose constraint (c) when $2p < q$ rather than only when $3p < q$ will become clear in the proof of Lemma \ref{4pquj}.)

Let $D$ be any cycle of $\mathcal{L}_2$. By Lemma \ref{lb}, we can choose a long bracket $u_1u_2u_3u_4$ in $D$. Consider the cycles $(u_1, w_{1, i})$---$(u_2, w_{2, i})$---$(u_3, w_{3, i})$---$(u_4, w_{4, i})$---$(u_1, w_{1, i})$ and $(u_1, w_{3, i})$---$(u_2, w_{4, i})$---$(u_3, w_{1, i})$---$(u_4, w_{2, i})$---$(u_1, w_{3, i})$ of $L$ on $A$. Both of them are rhombus switches for $G$ and the two pairs of cycles they connect belong to the two braids of $D$ in $G$.

Let $\mathcal{S}_\text{III}$ be the system of all such rhombus switches over all $i$ with $1 \le i < k$ and all cycles $D$ of $\mathcal{L}_2$. Then the switches of $\mathcal{S}_\text{III}$ are pairwise edge-disjoint.

Consequently, the switches in each one of our systems are pairwise edge-disjoint. This holds also for the union of all three systems: Two switches out of different systems cannot have any edges in common because their $y$-projections are edge-disjoint by virtue of constraints (a), (b), and (c).

Our setup is complete.

Flip all switches of $\mathcal{S}_\text{I}$ simultaneously. Since all of them are in separate subscarves, there will be no interference between them, and the net effect will be that within each of the first $\ell$ subscarves of $G$ some $p + q$ cycles will become stitched together into one single longer cycle. Denote the pseudotour of $L$ on $A$ thus obtained by $G^\star$.

We claim that the union $\mathcal{S}_\text{II} \cup \mathcal{S}_\text{III}$ is a connecting system of switches for $G^\star$.

Let us see why. Denote by $G^{\star\star}$ the union of $G^\star$ with all switches in the systems $\mathcal{S}_\text{II}$ and $\mathcal{S}_\text{III}$.

The switches of $\mathcal{S}_\text{II}$ ensure that, for each $i$, all edges of $G^\star$ in the same braid of the subscarf $G_i$ belong to the same connected component of $G^{\star\star}$.

Then the switches of $\mathcal{S}_\text{III}$ ensure that all edges of $G^\star$ in the same braid of the complete scarf $G$ belong to the same connected component of $G^{\star\star}$.

Finally, our initial flipping of all switches of $\mathcal{S}_\text{I}$ furthermore ensures that the entire union $G^{\star\star}$ is connected, as needed.

Therefore, $\mathcal{S}_\text{II} \cup \mathcal{S}_\text{III}$ is a connecting system for $G^\star$ which consists of pairwise edge-disjoint rhombus switches. By the rhombus switches lemma, we are done. \end{proof*}

\medskip

We go on to Theorem \ref{main}. We will need the following result of \cite{B}.

\medskip

\begin{theorem*}{\textbf{W}} (Supplement to Willcocks's conjecture.) Let $L$ be a skew free $(p, q)$-leaper. Then $L$ tours all boards both of whose sides are multiples of $2(p + q)$. \end{theorem*}

\medskip

The proof of the existence part of Willcocks's conjecture in \cite{B} is somewhat involved and we cannot reproduce it here. The gist is as follows: First we construct an explicit pseudotour of $L$ on the square board of side $2(p + q)$ out of twenty-eight pencils. (The definition of a pencil in a leaper graph is analogous to the definition of a pencil in a projection graph which we gave in Section \ref{proj-i}. It is the set of all edges joining corresponding cells in two rectangular subboards related by a suitable translation.) Then out of eight pencils we construct also an explicit connecting system of pairwise edge-disjoint rhombus switches for that pseudotour. Finally we apply the rhombus switches lemma. The main difficulty is to show that the system of switches does indeed connect the pseudotour.

The supplement above follows from the existence part of Willcocks's conjecture by the same strategy as in our discussion of Questions \textbf{A}---\textbf{D} in the introduction. This is also how we will establish Theorem \ref{main}.

We must stitch together some tours of the form described in the supplement to Willcocks's conjecture and some tours of the form described in Theorem \ref{4pq}. We have some degree of flexibility regarding the latter tours but not regarding the former ones. Thus we carry out the stitching in such a way that the interior structure of the former tours does not matter much.

Consider a board $A = I \times J$. Suppose, for concreteness, that $\min I = \min J = 0$. We call an edge of $L$ on $A$ a \emph{universal joint} when it has the form $(u, v)$---$(u - p, v + q)$ with $q - p \le u < q$ and $p \le v < 2p$.

\medskip

\begin{lemma} \label{uj} Suppose that there exists a tour of $L$ on the board of size $m \times n_\text{Joint}$ with a universal joint. Suppose, furthermore, that $L$ tours the board of size $m \times n$. Then $L$ also tours the board of size $m \times (n + n_\text{Joint})$. \end{lemma} 

\medskip

\begin{proof} Partition an integer interval $I$ of size $n + n_\text{Joint}$ into two subintervals $I'$ and $I''$ of sizes $n$ and $n_\text{Joint}$, respectively, and let $J$ be an integer interval of size $m$. Let $T'$ be a tour of $L$ on the board $A' = I' \times J$ and let $T''$ be a tour of $L$ with a universal joint $(u, v)$---$(u - p, v + q)$ on the board $A'' = I'' \times J$.

Observe that cell $(u - q, v - p)$ is only incident with two edges in the leaper graph of $L$ on $A'$, and so both of them are forced to be in $T'$. We are only left to flip the rhombus switch $(u, v)$---$(u - p, v + q)$---$(u - p - q, v - p + q)$---$(u - q, v - p)$---$(u, v)$ which connects $T'$ and $T''$. \end{proof}

\medskip

\begin{lemma} \label{4pquj} Theorem \ref{4pq} continues to hold when we furthermore require that the tour contains a universal joint. \end{lemma} 

\medskip

\begin{proof} For the knight, we verify this directly.

With all other skew free leapers, in the notation of the proof of Theorem \ref{4pq} we must ensure that: (a) The loom $\mathcal{L}_2$ contains an edge $e' = u$---$(u - p)$ with $q - p \le u < q$; (b) The pattern $P_1$ contains an edge $e'' = v$---$(v + q)$ with $p \le v < 2p$; and (c) The universal joint of $L$ on $A$ whose projections are $e'$ and $e''$ is not destroyed when we flip some switches as prescribed by the rest of the proof.

When $p = 1$, part (a) is immediate. Otherwise, when $p \ge 2$, since $p$ does not divide $q$ the integer interval $[q - p; q - 1]$ of size $p$ must contain vertices $u$ realising both parities of $\lfloor u/p \rfloor$. Thus it will also contain some vertex $u$ whose short edge in $\mathcal{L}_2$ lies on its left.

When $2p < q$, part (b) follows by constraint (c) in the proof of Theorem \ref{4pq}. Otherwise, when $2p > q$, it follows because the leftmost edges of the pattern $P_1$ form a translation copy of $\mathcal{B}_1$. (We introduced the graphs $\mathcal{B}_k$ in Section \ref{proj-ii}.)

For part (c), our universal joint cannot be in any switch of systems $\mathcal{S}_\text{I}$ and $\mathcal{S}_\text{II}$ by constraint (c) in the proof of Theorem \ref{4pq} when $2p < q$ and because all instances of $\mathcal{E}$ specified in constraints (a) and (b) of that proof are edge-disjoint from the aforementioned translation copy of $\mathcal{B}_1$ when $2p > q$. It cannot be in any switch of system $\mathcal{S}_\text{III}$, either, because its $x$-projection is short and the $x$-projections of all edges of these switches in $G$ are long. \end{proof}

\medskip

\begin{proof*}{Proof of Theorem \ref{main}} Let $m$ and $n$ be sufficiently large and even. Let $m'$ be the greatest multiple of $2(p + q)$ such that $m' \le m$ and $m' \equiv m \pmod{4pq}$, and define $n'$ similarly.

By the supplement to Willcocks's conjecture, there exists a tour of $L$ on the board of size $m' \times n'$. By Theorem \ref{4pq}, there exists a tour of $L$ on the board of size $m' \times 4pq$ as well. By Lemma \ref{4pquj}, we can safely assume that the latter tour contains a universal joint. Then by repeated application of Lemma \ref{uj} we obtain a tour of $L$ on the board of size $m' \times n$.

For the transition from the board of size $m' \times n$ to the board of size $m \times n$ we proceed analogously. \end{proof*}

\medskip

Let $\mu_\square$ be the smallest even positive integer such that $L$ tours the square board of side $n$ for all even positive integers $n$ with $n \ge \mu_\square$. Let also $\mu_\text{Even}$ be the smallest even positive integer such that $L$ tours all boards both of whose sides are even and greater than or equal to $\mu_\text{Even}$. Thus $\mu_\square \le \mu_\text{Even}$.

Theorem \ref{main} shows that both of $\mu_\square$ and $\mu_\text{Even}$ are well-defined for all $L$. Furthermore, from its proof we can extract explicit upper bounds for them.

Suppose, to begin with, that $L$ is distinct from the knight and the antelope, so that the extension $\mathcal{E}$ in the proof of Theorem \ref{4pq} is well-defined and of length $2q$.

Our observations following Theorem \ref{pt} show that we can construct projection graph tours of all sizes greater than or equal to $6q^2$ so that they admit splits of the same width as $\mathcal{E}$. By the proof of Lemma \ref{swb}, we see also that it suffices for the chain of extensions considered there to contain at least $\xi$ instances of $\mathcal{E}$, where $\xi$ is as in Section \ref{loom}.

Let $m_\text{I} = 6q^2 + (\xi + 2) \cdot 2q$. Then, in the notation of the proof of Theorem \ref{4pq}, so long as $m_i \ge m_\text{I}$ for all $i$ we can indeed construct patterns $P_i$ satisfying all three constraints (a), (b), and (c).

Let $m_\text{II} = 4\eta(m_\text{I} + 1)$, where $\eta$ is as in Section \ref{loom}, and suppose that $m \ge m_\text{II}$. Then we can ensure that $k$ has the same parity as $m$, $k \ge \ell$, and $m_i$ is odd with $m_i > m_\text{I}$ for all $i$.

Finally, let also $m_\text{III} = m_\text{II} + 4pq(p + q)$, and suppose that $m \ge m_\text{III}$ and $n \ge m_\text{III}$. Then in the proof of Theorem \ref{main} we obtain $m' > m_\text{II}$ and $n' > m_\text{II}$, and so Theorem \ref{4pq} will indeed apply as stated. Consequently, $\mu_\text{Even} \le m_\text{III}$.

The obvious upper bounds $\eta \le pq/2$ and $\xi \le pq/2$ from Section \ref{loom} are enough to check that $m_\text{III} < 5q^5$. We modify the calculations appropriately for the knight and the antelope in order to conclude that the same estimate for $m_\text{III}$ holds for them as well. Therefore, $\mu_\square \le \mu_\text{Even} < 5q^5$ for all skew free leapers $L$, and so both of $\mu_\square$ and $\mu_\text{Even}$ are $\mathcal{O}(q^5)$.

\section{Further Work} \label{further}

The tours of projection graphs are exciting objects of study in their own right, quite apart from their connection to leaper tours. So let us talk about them first.

Outside of the context of leapers, we must rethink what values of the parameters it makes sense to consider. Let $a$ and $b$ be positive integers with $a < b$. When $a$ and $b$ share a nontrivial common divisor, $\Pi(a, b, n)$ is never Hamiltonian. Thus suppose, from this point on throughout the rest of our discussion of projection graphs, that $a$ and $b$ are relatively prime.

When $a$ and $b$ are of opposite parities, they are also the parameters of a skew free leaper, and our analysis of Sections \ref{proj-i} and \ref{proj-ii} applies to them.

Otherwise, when both of $a$ and $b$ are odd, $\Pi(a, b, n)$ becomes bipartite. The two parts consist of its even and odd vertices, respectively. Thus it can only be Hamiltonian when $n$ is even.

Everything we said in Sections \ref{proj-i} and \ref{proj-ii} continues to apply here as well, with one exception: In the proofs of Lemmas \ref{pte} and \ref{ptd}, $a + b$ and $2b$ are not relatively prime anymore, but rather $\gcd(a + b, 2b) = 2$. Therefore, in this case Theorem \ref{pt} must be modified to state instead that $\Pi(a, b, n)$ is Hamiltonian for all sufficiently large even positive integers $n$.

Let us augment the definition of $\mu_\Pi$ so that when both of $a$ and $b$ are odd it becomes the smallest even positive integer such that $\Pi(a, b, n)$ is Hamiltonian for all even positive integers $n$ with $n \ge \mu_\Pi$. Then, as in Section \ref{proj-ii}, in all cases $\mu_\Pi < 6b^2 \in \mathcal{O}(b^2)$.

What is the exact value of $\mu_\Pi$? Or, alternatively, what are some more precise lower and upper bounds for it?

Then, more generally, what can we say about the positive integers $n$ such that $\Pi(a, b, n)$ is Hamiltonian?

For example, the proof of Lemma \ref{pst} shows that the smallest positive integer $n$ for which $\Pi(a, b, n)$ is Hamiltonian is $n = a + b$. Furthermore, Lemma \ref{pmt} shows that every multiple of $a + b$ yields a Hamiltonian projection graph, too. Thus it is natural to wonder what happens outside of the sequence of the multiples of $a + b$.

Let $\mu_\text{Div}$ be the smallest positive integer such that $a + b$ does not divide $\mu_\text{Div}$ and yet $\Pi(a, b, \mu_\text{Div})$ is Hamiltonian. 

By Lemma \ref{pkt}, if $a = 1$ and $b = 2$ then $\mu_\text{Div} = 4$.

\medskip

\begin{proposition} \label{pde} Suppose that $2a < b$. Then $\mu_\text{Div} = 3a + b$. \end{proposition} 

\medskip

\begin{proof} Let $I$ be an integer interval of size $n$.

Suppose first that $a + b < n < 3a + b$. Then there exist two vertices $u'$ in $\First(a, I)$ and $u''$ in $\Last(a, I)$ with $u'' - u' = a + b$. In every pseudotour of $\Pi(a, b, I)$, the forced edges of $u'$ and $u''$ form a cycle of length four.

Suppose, then, that $n = 3a + b$. Partition $I$ into subintervals $I_1$, $I_2$, $J$, $I_3$, and $I_4$ of sizes $a$, $a$, $b - a$, $a$, and $a$, respectively. Consider the union of the pencils $I_1 \circeq I_2$, $I_3 \circeq I_4$, $\First(b - 2a, J) \circeq \Last(b - 2a, J)$, and $\First(3a, I) \circeq \Last(3a, I)$. The argument that this union is a tour proceeds along the same lines as in the proofs of Lemmas \ref{sde} and \ref{sdo}. \end{proof}

\medskip

Thus if $2a < b$, then $a + b < \mu_\text{Div} < 2(a + b)$. The situation changes when $2a > b$.

\medskip

\begin{conjecture} \label{pdd} Suppose that $2a > b$. Let $d = b - a$ and $\alpha = \lfloor b/d \rfloor$. Then $\alpha(a + b) < \mu_\text{Div} \le \alpha(a + b) + d$ when $d \ge 2$ and $\mu_\text{Div} = (\alpha - 1)(a + b) + 1$ otherwise, when $d = 1$. \end{conjecture} 

\medskip

The author has verified Conjecture \ref{pdd} for all $107$ admissible pairs of parameters $a$ and $b$ with $d \le 7$ and $\alpha \le 7$.

We stumble upon more curious findings when we consider the number of tours that projection graphs possess.

Let $\mu_\text{Var}$ be the smallest positive integer such that $\Pi(a, b, \mu_\text{Var})$ admits at least two distinct tours. When no such positive integer exists, formally set $\mu_\text{Var} = \infty$. 

\medskip

\begin{conjecture} \label{pvt} The identity $\mu_\text{Var} = 2(a + b)$ holds for all pairs of parameters $a$ and $b$ with the exception of: (a) $a = 1$ and $b = 2$, when $\mu_\text{Var} = \infty$; (b) $a = 1$ and $b$ odd, when $\mu_\text{Var} = 2b$; and (c) $b - a = 1$ and $a \ge 2$, when $\mu_\text{Var} = 3(a + b)$. \end{conjecture} 

\medskip

The author has verified Conjecture \ref{pvt} for all $386$ pairs of parameters $a$ and $b$ with $a + b \le 50$.

Conversely, when $n$ is smaller than our conjectured value of $\mu_\text{Var}$ and $\Pi(a, b, n)$ is Hamiltonian, we should expect it to be uniquely so.

One corollary is that the tour of Proposition \ref{pde} should be unique in all cases except for $a = 1$ and $b = 3$. This is indeed true. To see why, consider an arbitrary pseudotour of $\Pi(a, b, 3a + b)$. In the notation of the proof of Proposition \ref{pde}, if our pseudotour omits any long edge of the pencil $I_2 \circeq I_3$, then forced edges will form a cycle of length six in it. Thus when $3a + b > 6$ a tour must necessarily contain this entire pencil, and the rest is straightforward.

Conjecture \ref{pvt} shows that we should expect to see some uniquely Hamiltonian projection graphs with relatively few vertices. Surprisingly, sometimes we observe unique tours in projection graphs of much greater sizes as well.

\medskip

\begin{conjecture} \label{put} Suppose that $2a > b$. Let $d = b - a$ and $\alpha = \lfloor b/d \rfloor$. Suppose, furthermore, that $\alpha$ is even and that $a \bmod d = b \bmod d > d/2$. Then the projection graph $\Pi(a, b, 2\alpha b - d)$ admits a unique tour. \end{conjecture} 

\medskip

The author has verified Conjecture \ref{put} for all $60$ admissible pairs of parameters $a$ and $b$ with $d \le 10$ and $\alpha \le 8$.

The tour of Conjecture \ref{put} is structured similarly to the tours we constructed in the proofs of Lemmas \ref{sde} and \ref{sdo}, as follows.

Let $r = a \bmod d = b \bmod d$ and $k = \alpha/2$. Let $I$ be an integer interval of size $2\alpha b - d$, and partition $I$ into subintervals $J'$, $D_1$, $D_2$, \ldots, $D_{4k - 1}$, $J''$ of sizes $2ka$, $d$, $d$, \ldots, $d$, $2ka$, respectively. Place one translation copy of $\mathcal{B}_k$ on vertices $J' \cup D_1 \cup D_2 \cdots \cup D_{2k}$ as well as one reflection copy of it on vertices $D_{2k} \cup D_{2k + 1} \cup \cdots \cup D_{4k - 1} \cup J''$. Then also add in the pencils $\First(d - r, D_j) \circeq \Last(d - r, D_{j + 2k})$ and $\Last(r, D_j) \circeq \First(r, D_{j + 2k})$ for all $j$ with $1 \le j \le 2k - 1$.

We can establish that this construction does indeed yield a tour by the same method as in the proofs of Lemmas \ref{sde} and \ref{sdo}.

One more engaging problem would be to determine exactly -- or at least to estimate somewhat more precisely -- the quantities $\eta$ and $\xi$ associated with looms, as defined in Section \ref{loom}.

With this, we return to leapers.

It would be very interesting to see a proof of the even area conjecture for boards with one odd and one even side.

Theorem \ref{4pq} already shows that $L$ tours some boards of this form. How close can we get to the full even area conjecture by combining together the tours currently at our disposal?

Fix an even positive integer $c$. We say that a board is \emph{plain} when either both of its sides are even or one of its sides is odd and the other one is a multiple of $c$. Then all boards which $L$ tours by the supplement to Willcocks's conjecture and Theorem \ref{4pq} are plain with $c = 4pq$.

\medskip

\begin{proposition} \label{plain} Suppose that $A$ admits a dissection into plain boards. Then $A$ itself is plain as well. \end{proposition}

\medskip

\begin{proof} We can assume without loss of generality that all pieces in the dissection are of sizes $2 \times 2$, $1 \times c$, and $c \times 1$.

Suppose that the height of $A$ is odd and let $A = I \times J$. Consider each subinterval of $I$ of size $c$ that is the $x$-projection of an odd number of pieces of size $1 \times c$ in the dissection. Every element of $I$ must belong to an odd number of such subintervals. That is only possible when $c$ divides the size of $I$, and so also the width of $A$. \end{proof}

\medskip

Proposition \ref{plain} with $c = 4pq$ shows that with the building blocks immediately available to us we could at best hope to construct tours of $L$ on boards with one odd side whose even side is a multiple of $4pq$. By Theorem \ref{4pq} and Lemmas \ref{uj} and \ref{4pquj}, $L$ does indeed tour all sufficiently large boards of this form.

We can do a tiny bit better if we replace $\mathcal{L}_2$ in Sections \ref{switch-ii} and \ref{leap} with an arbitrary loom $\mathcal{L}_k$ where $k \ge 2$. That would necessitate only minor changes in the argument, and it would also allow us to construct tours of $L$ on all sufficiently large boards with one odd side whose even side is a multiple of $2pq$.

But the complete even area conjecture still seems to lie a long way away.

Let $\mu_\text{Tour}$ be the smallest positive integer such that $L$ tours all boards of even area both of whose sides are greater than or equal to $\mu_\text{Tour}$. If no such positive integer exists, then formally set $\mu_\text{Tour} = \infty$. Thus $\mu_\square \le \mu_\text{Even} \le \mu_\text{Tour} + (\mu_\text{Tour} \bmod 2)$. Furthermore, $\mu_\text{Tour}$ is finite if and only if the even area conjecture holds for $L$.

To determine $\mu_\square$, $\mu_\text{Even}$, and $\mu_\text{Tour}$ exactly would probably be an intractable problem in the general case. However, we could still attempt to obtain some more precise lower and upper bounds for them.

By the closing remarks of Section \ref{leap}, both of $\mu_\square$ and $\mu_\text{Even}$ are $\mathcal{O}(q^5)$.

On the other hand, a corollary of one result of \cite{O} we cited in the introduction is that, when $L$ belongs to the exceptional family $q = p + 1$, it does not tour any board of height $2pq$. Thus for such leapers $2pq < \mu_\square$. Consequently, any upper bound for one of $\mu_\square$, $\mu_\text{Even}$, and $\mu_\text{Tour}$ of the form $\mathcal{O}(q^\vartheta)$ which applies to all skew free leapers must necessarily satisfy $\vartheta \ge 2$.

It is at least somewhat plausible that all three of $\mu_\square$, $\mu_\text{Even}$, and $\mu_\text{Tour}$ are in fact $\mathcal{O}(q^2)$.

One more natural question is as follows. Let $\mu_\text{Side}$ be the smallest positive integer such that $L$ tours some board of that height. (Or, equivalently, some board of that width.) What can we say about $\mu_\text{Side}$?

\medskip

\begin{conjecture} \label{side} The identity $\mu_\text{Side} = p + q$ holds for all skew free leapers $L$. \end{conjecture}

\medskip

The author has verified Conjecture \ref{side} for all $24$ skew free leapers $L$ with $p + q \le 15$.

We can resolve the analogous question for pseudotours right away. Let $\mu^\star_\text{Side}$ be the smallest positive integer such that there exists a pseudotour of $L$ on some board of that height.

\medskip

\begin{proposition} \label{ps} The identity $\mu^\star_\text{Side} = p + q$ holds for all skew free leapers $L$. \end{proposition} 

\medskip

\begin{proof} When $m < p + q$, every board of height $m$ contains at least one cell of degree at most one in the leaper graph of $L$. On the other hand, by Lemmas \ref{pst}, \ref{2pq}, and \ref{ls} there exists a scarf pseudotour of $L$ on the board of size $(p + q) \times 2pq$. \end{proof}

\medskip

Clearly, then also $\mu_\text{Side} \ge p + q$.

One promising approach to Conjecture \ref{side} would be to search for a scarf pseudotour of $L$ on a board of the form $(p + q) \times 2kpq$ which can be modified somehow so as to become a tour.

The work \cite{Kn} contains some more partial results relevant to Conjecture \ref{side}: Suppose that $L$ tours the board of size $(p + q) \times n$. Then if $L$ belongs to the exceptional family $p = 1$, necessarily $n \ge 4q + 2$. (Knuth conjectures also that in this special case the smallest positive integer $n$ for which a tour does exist is in fact $n = 6q$ when $q \ge 4$ and $n = 10$ when $q = 2$.) Furthermore, if $L$ belongs to the exceptional family $q = p + 1$, then for all $p \ge 4$ necessarily $n \ge q^2 + 3q - 2$. (With the lower bound raised up slightly when $p$ is even and when $p$ is a multiple of four.)

The other problems we have considered thus far admit pseudotour variants, too. For example, the pseudotour analogue of the even area conjecture does not appear to have been studied before. We can also define the quantities $\mu^\star_\square$, $\mu^\star_\text{Even}$, and $\mu^\star_\text{Tour}$ in the same way as $\mu_\square$, $\mu_\text{Even}$, and $\mu_\text{Tour}$ but with pseudotours instead of tours, and ask the exact same questions about them. Clearly, then $\mu^\star_\square \le \mu_\square$ and similarly for the other two pairs.

These pseudotour analogues are likely to be a lot more accessible, and work on them will perhaps shed some light on their more difficult tour counterparts as well.

The boards in the proof of Proposition \ref{ps} lead to some settings where scarves occur naturally.

\medskip

\begin{proposition} \label{use} Suppose that $L$ belongs to the exceptional family $p = 1$. Then there exists a unique pseudotour of $L$ on the board of size $(p + q) \times 2pq$. \end{proposition} 

\medskip

\begin{proof} The cells of the two outermost rows except for the four corner cells are of degree three in the leaper graph of $L$ and the other cells of the board are all of degree two. Therefore, for a pseudotour we must remove $2q - 2$ edges forming a perfect matching over the $4q - 4$ cells of degree three. \end{proof}

\medskip

\begin{proposition} \label{usd} Suppose that $L$ belongs to the exceptional family $q = p + 1$. Then there exists a unique pseudotour of $L$ on the board of size $(p + q) \times 2pq$. \end{proposition} 

\medskip

The proof is somewhat more involved this time around.

\medskip

\begin{proof} Let $A$ be a board of size $(p + q) \times 2pq$ and let $G$ be a pseudotour of $L$ on $A$. We puzzle out $G$ edge by edge.

Let us colour an edge of $L$ on $A$ green when we have deduced that it must belong to $G$ and red when we have deduced that it must be outside of $G$. We make all of our deductions according to the following pair of rules.

\smallskip

\emph{Rule \textbf{I}}. Suppose that we have already coloured all but two of the edges of some cell in red. Then we can colour its two remaining edges in green. 

\smallskip

\emph{Rule \textbf{O}}. Suppose that we have already coloured two of the edges of some cell in green. Then we can colour all of its remaining edges in red. 

\smallskip

We proceed step by step over the course of $2p$ consecutive steps numbered from $0$ to $2p - 1$. On each step $k$, we consider some subset of cells and we apply to each one of them Rule \textbf{I} when $k$ is even and Rule \textbf{O} when $k$ is odd.

These subsets are as follows. Label the columns of the board $0$, $1$, \ldots, $2pq - 1$ from left to right.

When $0 \le k < p$, on step $k$ we consider all cells in columns of the form $qk + i$ with $0 \le i < p - k$ as well as their reflections across the vertical axis of symmetry of the board. (Thus $2(p - k)(p + q)$ cells altogether.)

Otherwise, when $p \le k < 2p$, on step $k$ we consider all cells in columns of the form $qk - (p + q)i$ with $0 \le i \le k - p$ as well as their reflections across the vertical axis of symmetry of the board. (Thus $2(k - p + 1)(p + q)$ cells altogether.)

Following step $2p - 1$, every edge will be coloured either green or red, with the green edges forming the unique scarf pseudotour of $L$ on $A$. \end{proof}

\medskip

Our discussion of leaper tours is complete.

\bibliographystyle{amsplain}

\begin{aicauthors}
\begin{authorinfo}[nb]
  Nikolai Beluhov\\
  Bulgaria\\
  nikolai~[dot]~beluhov~[at]~gmail~[dot]~com
\end{authorinfo}
\end{aicauthors}

\end{document}